\documentclass[a4paper]{article}
\usepackage{graphicx}
\usepackage{amsmath,amssymb}
\usepackage[margin=20mm]{geometry}
\usepackage{url}
\newcommand{\diff}[1]{#1}
\title{Designing efficient interventions\\for pre-disease states using control theory}
\author{Makito Oku}
\date{}

\begin{document}
\maketitle

\begin{abstract} 
To extend healthy life expectancy in an aging society, it is crucial to prevent various diseases at pre-disease states. Although dynamical network biomarker theory has been developed for pre-disease detection, mathematical frameworks for pre-disease treatment have not been well established. Here I propose a control theory-based approach for pre-disease treatment, named Markov chain sparse control (MCSC), where time evolution of a probability distribution on a Markov chain is described as a discrete-time linear system. By designing a sparse controller, a few candidate states for intervention are identified. The validity of MCSC is demonstrated using numerical simulations and real-data analysis.
\end{abstract}

\section{Introduction}

Conventional medicine has mainly focused on patients who already have a particular disease. However, in an aging society \cite{muramatsu2011}, an increasing number of people are suffering from various diseases, such as cancer, dementia, and diabetes \cite{iijima2021}. Under such circumstances, treatment after onset is insufficient, and disease prevention is becoming increasingly important to extend healthy life expectancy and reduce medical costs. Although traditional preventive medicine has established general-purposed effective measures, such as improving lifestyle habits including exercise and diet, more targeted approaches are also needed to mitigate the rapid increase in disease incidence among older adults. The key to this is a deeper understanding of pre-disease states.

Pre-disease states are generally recognized as intermediate states between a healthy state and a disease state. Their ambiguous definition is often criticized as an obstacle to conducting rigorous research. Therefore, it is important to identify \textit{clinically significant} pre-disease states, where the risk of developing a certain disease is sufficiently high, and the benefits of intervention justify the costs. For example, metabolic syndrome is a clinically significant pre-disease state. It is not a disease state that reduces quality of life, but it has a high risk for diabetes, cardiovascular diseases, and so on. It can be diagnosed objectively and at low cost, using waist measurement, blood pressure measurement, and blood sampling. It can be restored to a healthy state at low cost through lifestyle changes without the use of drugs. The rationale for intervention at this stage is supported by reliable scientific evidence. Consequently, specific health checkups and specific health guidance have been established in Japan for the detection and treatment of metabolic syndrome, respectively \cite{tsushita2017}.

Previous studies used mathematical approaches to find other candidates for clinically significant pre-disease states. For example, dynamical network biomarker (DNB) theory has been developed for pre-disease detection, focusing on fluctuations in biological systems \cite{chen2012, aihara2022}. DNB theory models the dynamics of a complex biological system as a multivariate noisy dynamical system \cite{oku2018}, and assumes that healthy and disease states correspond to two distinct equilibrium points of the dynamical system, which is often explained using a double-well potential model as an analogy \cite{chen2012, aihara2022}. A representative scenario is that the initial state is the healthy state, and a bifurcation parameter changes gradually. As a saddle-node bifurcation approaches, the system state changes to a pre-disease state. At the pre-disease state, noise tolerance (also called resilience) declines, and fluctuations become larger, which is called an early warning signal (EWS) \cite{scheffer2009, scheffer2012}. After the bifurcation, the system state inevitably moves to another equilibrium point corresponding to the disease state. DNB theory provides a practical data analysis method to capture EWS from multivariate real data, enabling the identification of pre-disease states \cite{koizumi2019, yonezawa2024} and associated variables \cite{akagi2025}.

There are many other mathematical methods that may be suitable for pre-disease detection. First of all, machine learning techniques are undoubtedly useful for onset prediction \cite{nanri2015, segan2023}, patient stratification \cite{Jeong2024, yoshimura2025}, early diagnosis \cite{fujiki2025, kitajima2025}, and obtaining objective biomarkers \cite{zhu2022, sawai2025}. It should be noted that onset prediction and patient stratification are essentially the same as pre-disease detection in that they identify high-risk individuals. Next, there are many other EWS-based methods besides DNB theory \cite{dakos2008, veraart2011, dakos2012, oku2023, masuda2024}. Moreover, energy landscape analysis (ELA) is a multivariate time-series data analysis method, which can reveal state transitions among multiple stable states using an Ising model \cite{watanabe2014, ezaki2017, masuda2025}. A recent study reported that, by analyzing health checkup data using ELA, a novel pre-disease state prior to diabetes onset was found \cite{ito2025}. Finally, Reservoir computing is a mathematical model that enables various tasks, such as time-series forecasting and classification, with a much lower learning cost than deep learning models \cite{Jaeger2002, maass2002, tanaka2022, matsumoto2023}. A recent study developed an anomaly detection method using reservoir computing \cite{tamura2025}, which will be useful for detecting irregularities in nearly periodic dynamics, such as arrhythmia and locomotive syndrome. 

However, previous studies primarily focused on pre-disease detection, and mathematical frameworks for pre-disease treatment have not been well established. Indeed, control theory has been integrated into DNB theory (hereafter referred to as DNB intervention theory) to provide a systematic strategy to restrict candidate variables for interventions \cite{yasukata2023, shen2023, shen2025}. As a proof of concept, 147 DNB genes obtained from a mouse model of metabolic syndrome \cite{koizumi2019} were ranked by the DNB intervention theory \cite{yasukata2023, shen2023, shen2025}, and top 10 genes were selected. After excluding two non-coding genes, 8 genes were experimentally investigated \cite{akagi2025}. Among them, the third-ranked gene named \textit{Vasa} in flies and \textit{DDX4} in humans was identified as a novel gene involved in metabolic processes \cite{akagi2025}. This indicates that the DNB intervention theory is effective in screening therapeutic targets for pre-disease states. However, it is currently limited to drug therapy, and other non-drug therapies, such as exercise therapy and diet therapy, are not covered.

Another study on mathematical pre-disease treatment is an optimal path planning framework using a Bayesian model \cite{nakamura2021}. This framework suggests an optimal order of improvement for each variable, such as body weight and blood pressure. For example, a person may be recommended to first reduce body weight without worrying about blood pressure, and then be advised to lower blood pressure by taking antihypertensive drugs and suppressing salt intake. Although this framework can be applied to non-drug therapies, such an extreme approach may not be acceptable in clinical practice, even if it is theoretically optimal. Therefore, it is necessary to develop more practical pre-disease treatment methods based on mathematical theories.

On the other hand, reinforcement learning (RL) \cite{sutton2018} has been widely used to obtain optimal behavioral strategies (called policies) in diverse tasks, from game playing \cite{mnih2015} to autonomous driving \cite{elallid2022} and autonomous robots \cite{kaufmann2023, Radosavovic2024}. It is also utilized for artificial intelligence alignment \cite{gabriel2020}, which is called RL from human feedback \cite{bai2022}. One of the key strengths of RL is its ability to operate in complex systems with unknown dynamics. Through interaction with the environment, an agent iteratively updates its policy to maximize the expected cumulative reward. The flexibility and versatility of RL comes from its mathematical foundation, known as Markov decision process (MDP) model, which describes stochastic state transitions depending on both the current state and action. MDP allows for the derivation of phenomenological models in a data-driven manner, even when the underlying mechanistic models are difficult to obtain.

RL is expected to be suitable for designing efficient intervention strategies for pre-disease states. Indeed, many studies have used RL to implement healthcare systems for personalized and adaptive interventions \cite{yom-tov2017, gonul2021, trella2022, deliu2024}. However, it is usually difficult to test various treatments, such as experimental drug administration and unestablished dietary restrictions, on people who are not sick. This prevents a sufficient exploration of the policy space. In contrast, there are many time-series data of personal health records, which can be obtained from regular health checkups, wearable devices, and so on. Therefore, such state-only data may help construct a suboptimal policy for pre-disease treatment. Mathematically, actionless MDP model corresponds to Markov chain model.

In this paper, I propose a control theory-based approach for pre-disease treatment, named Markov chain sparse control (MCSC). Briefly, it consists of the following steps. First, given time-series data of the health state, the state is discretized. Next, the transition matrix of a Markov chain is calculated using the discretized data. A key point here is to describe the time evolution of the probability distribution on the Markov chain as a discrete-time linear system. In other words, we do not consider individual realizations or trajectories of the stochastic process, as in the Langevin equation, but rather focus on how the probability distribution changes over time, as in the Fokker-Planck equation. This technique allows us to obtain a linear and deterministic model, even when the original dynamical system is nonlinear and stochastic. Finally, a controller is designed with sparse regularization to identify a small number of candidate states for intervention. This allows for targeted interventions to prevent the disease onset.

The rest of the paper is organized as follows. Section 2 describes MCSC in more detail. Section 3 describes numerical simulation settings and real data analyzed in this study. Section 4 presents the results of the numerical simulations and real data analysis. Section 5 provides a discussion. Section 6 concludes the paper.

\section{MCSC}

\subsection{Data assumptions}

\begin{figure}[t]
\centering
\includegraphics[width=0.8\hsize]{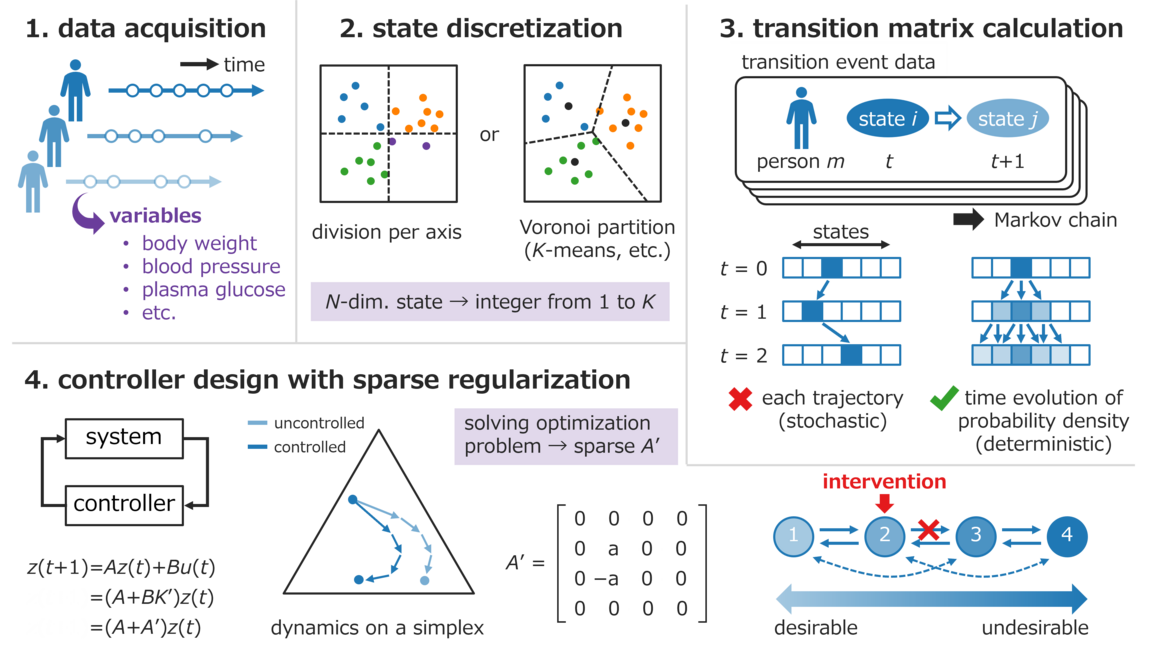}
\caption{Overview of MCSC for pre-disease treatment.}
\label{fig01}
\end{figure}

Figure~\ref{fig01} shows an overview of MCSC for pre-disease treatment. MCSC is developed for analyzing various time-series data of personal health records, such as regular health checkup data. Let $N$, $M$, and $T$ denote the number of observable variables, individuals, and time points, respectively. Variables are basically assumed to be numerical, but categorical variables can be included by assigning a unique numerical value to each category. Although $N$, $M$, and $T$ are usually greater than 1, $N$ and $M$ can be 1, but $T$ must be greater than 1. Missing values are allowed, but caution is needed when there is a significant imbalance in the number of data points between individuals. It is better to exclude variables with many missing values. Observations are assumed to be taken at approximately regular time intervals. Otherwise, resampling should be considered.

Time series of repeated measures are basically assumed, that is, the same individuals are repeatedly measured at multiple time points. However, MCSC can also be applied to time series of independent measures, where different individuals are measured at different time points. In that case, matching of individuals at different time points is needed, which will be explained in Section \ref{sec_matching}.

\subsection{State discretization}

Let $x(m,t)=(x_1(m,t), x_2(m,t), \ldots, x_N(m,t))^{\top}\in\mathbb{R}^N$ denote the health state of individual $m$ at time $t$. Instead of considering a dynamical system $x(m,t+1)=f(x(m,t))$ with $f:\mathbb{R}^N\to\mathbb{R}^N$, the state space is divided into $K$ disjoint subsets as follows:
\begin{equation}
\mathbb{R}^N = \bigcup_{k=1}^K S_k,\quad S_k\cap S_l=\emptyset\ (\forall k,\ l,\ k\neq l).
\end{equation}

Each state $x$ is mapped to an integer $k$ that satisfies $x\in S_k$. By denoting the map from $x$ to $k$ as $g:x(m,t)\in\mathbb{R}^N\to\tilde x(m,t)\in\{1,\ldots,K\}$, the discretized data are expressed as $\{g(x(m,t))\}=\{\tilde x(m,t)\}$.

In principle, any division is possible, but for practical reasons, $K$ should not be too large, and $\tilde x$ should not be overly concentrated at a single value.

\subsubsection{Division per axis}

A simple approach is to divide the domain of each variable into two or more subregions, and use their Cartesian product. For example, in ELA, the value of each variable is converted to a binary value in $\{0,1\}$ (or $\{-1,1\}$) depending on whether it is above or below a threshold (for example, the mean of the variable) \cite{watanabe2014, ezaki2017, masuda2025, ito2025}. After that, each combination of the binary values (for example, $000111_2=7$) is used to represent each subset $S_k\subset\mathbb{R}^N$. Because $K=2^N$, a small number of variables should be selected in advance for ELA. The same approach can be used for the present method. Moreover, although ternary ELA or $n$-ary ELA have not been developed, each axis can \diff{be} divided into three or more subregions when $N$ is sufficiently small (for example, $N\leq 4$).

\subsubsection{Voronoi partition}

More efficient space partitioning methods can be used. It is reasonable to allocate larger-volume subsets to low-density areas while using smaller-volume subsets for high-density areas. To achieve this, $K$ representative points $\{x^{(1)}, \ldots, x^{(K)}\}$ are chosen from the data distribution, and each state $x$ is assigned to the nearest representative point as $g(x)=\arg\min_k d(x,x^{(k)})$, where $d$ denotes a distance function. In this case, each subset is defined as $S_k=\{x\mid g(x)=k\}$, and the partition corresponds to a multivariate Voronoi diagram induced by the representative points.

There are several methods for selecting the representative points. The simplest way is random sampling from the observed data points without replacement. Another way is to use $K$-means clustering. Each center of the $K$ clusters is taken as a representative point. Note that the representative points do not have to be exactly the same as one of the observed data. Furthermore, modern dimension reduction techniques, such as t-distributed stochastic neighbor embedding (t-SNE) \cite{Maaten2008} and uniform manifold approximation and projection (UMAP) \cite{mcinnes2018}, can be used to convert the original data to low-dimensional data. Applying $K$-means to the low-dimensional data may help select more meaningful representative points. The assignment of each state to the closest representative point is also performed in the low-dimensional space.

Note that $K$-means here is not used to find natural clusters. The resulting clusters are not necessarily contiguous, and the cluster boundaries may pass through high-density areas. Rather than facilitating the interpretation of data, $K$-means here is used only to discretize the state space to an arbitrary resolution, as in quantization in analog-to-digital conversion.

\subsubsection{Other clustering methods}

Other clustering methods, such as hierarchical clustering, can also be used to define more flexible boundaries between subsets. However, many clustering methods are categorized as unsupervised machine learning and may not have an established way to assign a cluster number $\tilde x$ to each new data point $x$. One solution is to search for the nearest point in the training data to the input data point and use its cluster number. Another solution is to construct a classifier, such as a support vector machine, to predict the cluster number for the input data point.

\subsection{Calculation of the transition matrix of a Markov chain}

Next, the transition matrix $A$ of a Markov chain is calculated using the discretized data. In the rest of this paper, we simply refer to discretized states as “states.” $A_{ji}$ represents the transition probability from state $i$ to state $j$ in one step, that is, $P(\tilde x(m,t+1)=j\mid \tilde x(m,t)=i)$. The transition matrix $A$ is a $K\times K$ nonnegative matrix, and each column of $A$ must sum to 1, making it an element of the $(K-1)$-dimensional standard simplex $\Delta^{K-1}=\{p=(p_1,\ldots,p_K)^{\top}\mid \sum_ip_i=1,\ p_i\geq 0,\ \forall i\}$.

Let $q=(m,t,i,j)$ denote an event where individual $m$ is in state $i$ at time $t$ and changes to state $j$ at time $t+1$. The case of $i=j$ is allowed as self-transition. The set of all events $Q$ can be written as follows:
\begin{equation}
Q=\{q=(m,t,i,j)\mid \tilde x(m,t)=i,\ \tilde x(m,t+1)=j\}.
\end{equation}

\subsubsection{Relative frequency}

There are several methods to estimate $A$ from the event data $Q$. The simplest way is to use the relative frequency as follows:
\begin{equation}
A_{ji}=\frac{|\{q\mid q_3=i,\ q_4=j,\ q\in Q\}|}{|\{q\mid q_3=i,\ q\in Q\}|},\quad\forall i,j.
\end{equation}

\subsubsection{Damping factor}

If the Markov chain is not ergodic, it is reasonable to introduce a damping factor $\alpha\in(0,1)$ as in the PageRank algorithm as follows:
\begin{equation}
A\leftarrow \alpha\,A + \frac{1-\alpha}{K}\,O,
\end{equation}
where $O\in\mathbb{R}^{K\times K}$ is a matrix of ones. By rewriting $\alpha=1/(1+\varepsilon\,K)$, an alternative expression is obtained as follows:
\begin{equation}
A\leftarrow \frac{1}{1+\varepsilon\,K}\,(A + \varepsilon\,O).
\end{equation}

\subsubsection{Weighting}

If the number of data points varies widely between individuals, it may be better to assign smaller weights to individuals with more data points and vice versa, as follows:
\begin{align}
w(m,i)&=\frac{1}{|\{q\mid q_1=m,\ q_3=i,\ q\in Q\}|+\eta},\\
A_{ji}&=\frac{\sum_m w(m,i)\,|\{q\mid q_1=m,\ q_3=i,\ q_4=j,\ q\in Q\}|}{\sum_m w(m,i)\,|\{q\mid q_1=m,\ q_3=i,\ q\in Q\}|},\quad\forall i,j,
\end{align}
where $\eta>0$ is a regularization parameter.

\subsubsection{Smoothing}

If the estimated transition matrix $A$ is noisy and potentially unreliable, smoothing should be considered. For example, close states are often expected to have similar transition profiles. If that is the case, the distance matrix $D\in\mathbb{R}^{K\times K}$ with $D_{ij}=d(x^{(i)},x^{(j)})$ can be used to perform any distance-based smoothing, such as Nadaraya-Watson kernel regression. If a radial basis function with a parameter $\gamma>0$ is used as a kernel function, the $i$-th column of matrix $A$, denoted as $a^{(i)}=(a_{1i}, \ldots, a_{Ki})^{\top}\in\Delta^{K-1}$, becomes as follows:
\begin{equation}
a^{(i)}\leftarrow\frac{\sum_j\exp(-\gamma\,D_{ij}^2)\,a^{(j)}}{\sum_j\exp(-\gamma \,D_{ij}^2)}.
\end{equation}

\subsubsection{Resetting}

Time series data of personal health records are often non-stationary due to aging. Each person’s health status tends to deteriorate over time. When such directionality in state transitions is strong, adding a loop from the final state back to the initial state can help make the Markov chain stationary. In other words, periodic boundary conditions are applied on the time axis. Let $t_{\min}^{(m)}$, $t_{\max}^{(m)}$ denote the first and last time points of observations for person $m$, respectively. Resetting the state for each person corresponds to substituting the event set $Q$ as follows:
\begin{equation}
Q\leftarrow Q\cup\left\{\left(m ,t_{\max}^{(m)}, \tilde x\left(m, t_{\max}^{(m)}\right), \tilde x\left(m, t_{\min}^{(m)}\right)\right)\right\}.
\end{equation}

The destination after the reset can be made independent of individuals by adding a dummy state $K+1$. The final state of each person is connected to the dummy state, and the dummy state is connected to the initial state of each person.

\subsubsection{Time evolution of probability distribution}

After $A$ is obtained, let us consider the time evolution of the probability distribution of the Markov chain. \diff{An important point is that even though each trajectory of the Markov chain is stochastic, the dynamics of the probability distribution becomes deterministic.} Let $z(t)=(z_1(t), \ldots, z_K(t))^{\top}\in\Delta^{K-1}$ denote the probability distribution on the Markov chain at time $t$. The time evolution of $z(t)$ can be described by the following discrete-time linear system:
\begin{equation}
z(t+1)=A\,z(t).
\end{equation}

If the Markov chain is ergodic (i.e., irreducible and aperiodic), $z(t)\to z^*$ as $t\to\infty$, where the unique stationary distribution $z^*$ is the eigenvector of $A$ corresponding to the eigenvalue 1. This equation can be rewritten in a $(K-1)$-dimensional system as described in Appendix~\ref{seca2}.

\subsubsection{Matching}\label{sec_matching}

When MCSC is applied to time series of independent measures, such as single cell RNA-sequencing (scRNA-seq) data taken from multiple time points, we need to match individuals (persons, cells, and so on) at different time points. The problem to be solved is to obtain a transition matrix $A$ given the time series of the probability distribution $z(1),\ldots,z(T)$ and the distance matrix $D$. Since $A$ must be nonnegative, a simple estimation using the pseudo-inverse matrix $A=Z_+Z_-^{\dagger}$ is inappropriate, where $Z_+=[z(2)\cdots z(T)]$, $Z_-=[z(1)\cdots z(T-1)]$, and $ X^{\dagger}=(X^{\top}X)^{-1}X^{\top}$.

A practical solution is to use optimal transport (OT), which is also known as earth mover distance \cite{rubner2000} or Wasserstein distance, under the assumption that $A$ is time-invariant. Let $F(t)\in\mathbb{R}^{K\times K}$ denote the OT plan from $z(t)$ to $z(t+1)$ with respect to $D$, where $F_{ij}(t)$ denotes the flow from $i$ to $j$. Then, the average transport plan $\bar F$ is calculated as follows:
\begin{equation}
\bar F=\frac{1}{T-1}\sum_{t=1}^{T-1}F(t).
\end{equation}

After that, each column of $\bar F^{\top}$ is divided by its summation, resulting in the transition matrix $A$. If the result of the free-run simulation of Eq. (10) is largely inconsistent with the observations, it is better to apply smoothing to $z(1),\ldots,z(T)$ in the time direction in advance.

\subsection{Designing a controller with sparse regularization}

In control theory, we usually aim to stabilize a system by adding a feedback loop. Although an ergodic Markov chain is stable by itself, our goal is to move its stationary distribution to a more desirable position within the standard simplex $\Delta^{K-1}$. The state space representation of the system with feedback control is described as follows:
\begin{equation}
z(t+1)=A\,z(t)+B\,u(t),
\end{equation}
where $B\in\mathbb{R}^{\diff{K\times L}}$ is the input matrix, and $u(t)\in\mathbb{R}^{\diff{L}}$ is the control input at time $t$. \diff{$L$ is the dimension of the control input. Given that $z(t)$ represents the probability distribution of the Markov chain, the control input $u(t)$ is also real-valued.} $z(t+1)$ must be in $\Delta^{K-1}$. This equation should be interpreted as a virtual model because the probability distribution $z$ cannot be observed and manipulated \diff{by each medical doctor} in reality.

To make the system closed, the gain matrix $K'\in\mathbb{R}^{\diff{L \times K}}$ is introduced. Then, the system can be rewritten as follows:
\begin{align}
u(t)&=K'\,z(t),\\
z(t+1)&=(A+BK')\,z(t),\\
&=(A+A')\,z(t),
\end{align}
where $A'=BK'\in\mathbb{R}^{K\times K}$ with $\sum_iA'_{ij}=0\ (\forall j)$ and $A'_{ij}+A_{ij}\geq 0\ (\forall i,j)$. \diff{Although $A'$ can be defined directly without using $B$ and $K'$, it would be better to interpret $A'$ as $BK'$ for clarifying the relationship with control theory.} For practical reasons, $A'$ should be as sparse as possible. If $A'_{ji}>0$, the transition from state $i$ to state $j$ is promoted, and if $A'_{ji}<0$, the transition is suppressed. Thus, while the model itself is virtual, the controller derived from it provides practical insight into which state to intervene in and how. For example, consider a medical doctor acting as a controller. The doctor monitors the condition of a patient during outpatient visits. If $\diff{A'_{ji}<0}$ and the patient enters state $i$, the doctor may recommend lifestyle changes or prescribe medication to prevent the patient from transitioning to state $j$. Because MCSC is based on a Markov chain model rather than an MDP model, specific interventions must be considered on a case-by-case basis depending on the situation.

\subsubsection{\diff{Objective function}}

Let $r=(r_1,\ldots,r_K)^{\top}\in\mathbb{R}^K$ denote a reward vector, $z'\in\Delta^{K-1}$ denote the stationary distribution of $A+A'$, and $\lambda_1,\lambda_2\geq 0$ denote regularization parameters. The controller then is designed to maximize the following objective function $G$:
\begin{equation}
G=r^{\top}z' -\lambda_1\|\mathrm{vec}(A')\|_0 -\lambda_2 \sum_{i,j,A_{ij}\neq 0}\left|\log\frac{A_{ij}+A'_{ij}}{A_{ij}}\right|,
\end{equation}
where $\|\mathrm{vec}(A')\|_0$ is the $L_0$ norm of the vectorized matrix $A'$. The first term is the expected reward. The second term represents a sparse regularization. The third term is the quantitative cost of modifying the transition matrix measured by log fold-changes. Although many variations can be considered, this study investigated the behavior of $G$ in the above form as a first step.

If one attempts to control the probability distribution at a finite horizon $t=\tau$ rather than the infinite horizon $t\to\infty$, $z'$ should be redefined as $z'=(A+A')^{\tau-1}z(1)$.

\subsubsection{\diff{Optimization method}}

In this study, a greedy algorithm was used to optimize the objective function $G$. For simplicity, only suppression of off-diagonal elements of $A$ was considered, and the reduced probability was added to the self-transition of the source state. To reduce computational costs, only the top 20~\% off-diagonal elements were included in the intervention candidate set. Next, each candidate was independently suppressed slightly, and those that showed a decrease in the objective function value were removed from the candidate set. After that, an exhaustive search was performed to identify the best intervention site and the amount of suppression, and $A'$ was updated accordingly. By repeating this, $A'$ was iteratively updated until no improvement in $G$ was observed. The amount of suppression was chosen from a predefined set $H=\{h_i\}$. For example, if $H=\{0.5, 0.8, 0.9\}$, suppression levels of 50~\%, 80~\%, and 90~\% were tested. \diff{They should be in a practically acceptable range. For example, 99~\% suppression may not be achievable in most cases because strong medications cannot be used for pre-disease states and some individuals may not adhere to the recommendations for lifestyle changes.}

\section{Numerical simulation settings and real data}

\subsection{One-dimensional double-well model}

The following one-dimensional double-well model was used:
\begin{align}
U(x)&=x^4-2x^2,\\
f(x)&=-\frac{\text{d}U(x)}{\text{d}t}=-4x^3+4x,\\
\text{d}x&=f(x)\text{d}t+\sigma\,\text{d}W,
\end{align}
where $U$, $f$, $\sigma$, and $W$ are a potential function, a flow function, noise intensity, and a Wiener process, respectively. The stochastic differential equation (SDE) was numerically solved by using the Euler-Maruyama method as follows:
\begin{equation}
x(t+\Delta t)=f(x(t))\Delta t+\sqrt{\Delta t}\,\xi(t),
\end{equation}
where $\Delta t$ is the time interval, and $\xi(t)$ is a random variable following the standard normal distribution.

The default parameters for the one-dimensional double-well model were as follows: $K=20$, $T=10^4$, $\Delta t=0.01$, $\sigma=1$, $\varepsilon=10^{-10}$, and $H=\{0.5, 0.8, 0.9\}$. The regularization parameters were varied as $(\lambda_1,\lambda_2)\in\{(0.1,0.1),(0.05,0.05),(0.01,0.01)\}$. For checking the smoothing effect, $T$ was reduced to $10^3$, and $\gamma$ was set to 20 or 5. \diff{$r_i=1$ for $i\leq 10$, and $r_i=-1$ otherwise.}

\subsection{Two-dimensional double-well model}

The following two-dimensional double-well model was used:
\begin{align}
U(v)&=\frac{1}{4}(x+y)^4-(x+y)^2+2(x-y)^2,\\
f(v)&=-\nabla U(v)=\left(\begin{array}{c}-(x+y)^3-2x+6y\\-(x+y)^3-2y+6x\end{array}\right),\\
\text{d}v&=f(v)\text{d}t+\sigma I\text{d}W,
\end{align}
where $v=(x,y)^{\top}$. The potential function $U$ was obtained from $U'(x',y')=x'^4-2x'^2+4y'^2$ and variable transformations of $x=(x'+y')/\sqrt{2}$ and $y=(x'-y')/\sqrt{2}$. The SDE was solved by using the Euler-Maruyama method.

The default parameters for the two-dimensional double-well model were as follows: $T=10^4$, $\Delta t=0.05$, $\sigma=1$, $\varepsilon=10^{-10}$, $\lambda_1=0.005$, $\lambda_2=0.005$, and $H=\{0.5, 0.8, 0.9\}$. The number of states $K$ was 100 when each axis was divided into 10 subregions, and it was set to 50 for Voronoi partitions. \diff{For checking the effect of $K$, $K$ was varied as $K\in\{20,40,60,80,100\}$, and both $\lambda_1$ and $\lambda_2$ were changed to 0.003 because otherwise no intervention was suggested for $K=100$. $r_i=1$ for $x\leq 0\land y\leq 0$, $r_i=-1$ for $x>0\land y>0$, and $r_i=0$ otherwise.}

\subsection{Two-dimensional branching-flow model}

The following two-dimensional branching-flow model was used as a toy model to imitate the Waddington landscape \cite{waddington1957, wang2011}:
\begin{align}
U(v)&=h_1(y)\cos(2\pi x)+h_2(y)\cos(4\pi x)+h_3(y)\cos(8\pi x)+\theta y,\\
h_1(y)&=\max(2y-1,0),\\
h_2(y)&=1-|2y-1|,\\
h_3(y)&=\max(1-2y,0),\\
f(v)&=-\nabla U(v)=\left(\begin{array}{c} h_1(y)2\pi\sin(2\pi x)+h_2(y)4\pi\sin(4\pi x)+h_3(y)8\pi\sin(8\pi x)\\-\theta\end{array}\right),\\
\text{d}v&=f(v)\text{d}t+\left[\begin{array}{cc}\sigma&0\\0&0\end{array}\right]\text{d}W,
\end{align}
where $v=(x,y)^{\top}$. $h_1(y)$, $h_2(y)$, and $h_3(y)$ are membership functions of fuzzy sets, which satisfy $h_1(y)+h_2(y)+h_3(y)=1\ (\forall y\in[0,1])$. The SDE was solved by using the Euler-Maruyama method. The initial position for each trial was $v(0)=(0.5,1)^{\top}$.

The parameters for the two-dimensional branching-flow model were as follows: $K=80$, $M=100$, $T=100$, $\Delta t=0.001$, $\sigma=0.8$, $\varepsilon=10^{-10}$, $\lambda_1=0.0005$, $\lambda_2=0.0005$, and $H=\{0.5, 0.8, 0.9\}$. The parameter $\theta$ was set to $1/((T-1)\Delta t)$ so that each trajectory reached the $y=0$ line at 100 steps. \diff{$r_i=1$ for the discretized state closest to the target state, and $r_i=-1$ otherwise.}

\subsection{Three-dimensional chaotic attractors}

As examples of three-dimensional chaotic attractors, Lorenz attractor \cite{lorenz1963} and Rössler attractor \cite{rossler1976} were used. The Lorenz attractor is described as follows:
\begin{align}
\frac{\text{d}x}{\text{d}t}&=\sigma(y-x),\\
\frac{\text{d}y}{\text{d}t}&=x(\rho-z)-y,\\
\frac{\text{d}z}{\text{d}t}&=xy-\beta z,
\end{align}
where $\rho$, $\sigma$, and $\beta$ are the parameters of the Lorenz attractor. Standard parameter values were used in this study: $\rho=28$, $\sigma=10$, and $\beta=8/3$.

The Rössler attractor is described as follows:
\begin{align}
\frac{\text{d}x}{\text{d}t}&=-y-z,\\
\frac{\text{d}y}{\text{d}t}&=x+ay,\\
\frac{\text{d}z}{\text{d}t}&=b+z(x-c),
\end{align}
where $a$, $b$, and $c$ are the parameters of the Rössler attractor. Standard parameter values were used in this study: $a=0.1$, $b=0.1$, and $c=14$.

To solve the ordinary differential equations, the fourth-order Runge-Kutta method was used as follows:
\begin{align}
v(t+\Delta t)&=v(t)+\frac{k_1+2k_2+2k_3+k_4}{6},\\
k_1&=f(v(t)),\\
k_2&=f\left(v(t)+\frac{\Delta t}{2}\,k_1\right),\\
k_3&= f\left(v(t)+\frac{\Delta t}{2}\,k_2\right),\\
k_4&= f(v(t)+\Delta t\,k_3),
\end{align}
where $v=(x,y,z)^{\top}$ and $f(v)=\text{d}v/\text{d}t$. Initial transient dynamics for 1000 steps was removed.

The other parameters for three-dimensional chaotic attractor models were as follows: $K=50$, $T=20,000$, $\Delta t=0.01$, $\varepsilon=10^{-10}$, and $H=\{0.5\}$. Here, 0.8 and 0.9 were excluded from $H$ because they resulted in meaningless solutions where the orbit was restricted to a single point. The regularization parameters were set as $\lambda_1=\lambda_2=0.01$ for the Lorenz attractor and $\lambda_1=\lambda_2=0.002$ for the Rössler attractor. \diff{For the Lorenz attractor, $r_i=1$ for $x>0$, and $r_i=-1$ otherwise. For the Rössler attractor, $r_i=1$ for $z<1$, and $r_i=-1$ otherwise.}

\subsection{Public health checkup data}

National Database (NDB) open data are publicly available, provided by the Ministry of Health, Labour and Welfare, Japan. These are aggregated data of medical receipts and specific health checkups. It is not possible to identify individuals. In this study, the specific health checkup data of fiscal year 2021 from the ninth NDB open data were analyzed.

The number of areas (secondary medical care areas) was 335. For example, Tokyo was divided into 13 areas. There were 16 variables: body mass index (BMI), waist circumference (WC), fasting plasma glucose (FPG), hemoglobin A1c (HbA1c), systolic blood pressure (SBP), diastolic blood pressure (DBP), triglycerides (TG), high-density lipoprotein (HDL), low-density lipoprotein (LDL), aspartate aminotransferase (AST), alanine aminotransferase (ALT), $\gamma$-glutamyl transpeptidase ($\gamma$-GTP), hemoglobin, non-fasting blood pressure, creatinine, and estimated glomerular filtration rate. The last four variables contained missing values denoted as 0, and were excluded from analysis. There were seven age groups: 40--44, 45--49, 50--54, 55-59, 60--64, 65--69, and 70--74 years old. Data for males and females were provided \diff{separately}.

Although these data were not individual-based, they can be viewed as time series of average health status for men and women in each area. Under this interpretation, the number of virtual individuals was $M=670$, and the number of virtual time points was $T=7$. Because \diff{principal component analysis (PCA)} revealed large differences between males and females, they were analyzed separately, that is, $M=335$ for male data and $M=335$ for female data.

The other parameters for the public health checkup data were as follows: $K=10$, $\varepsilon=10^{-10}$, $\lambda_1=0.001$, $\lambda_2=0.001$, and $H=\{0.5\}$. Here, 0.8 and 0.9 were excluded from $H$ because it is almost impossible in practice to suppress a particular health state transition by more than 50~\%. \diff{$r_i=-1$ for $i=10$, and $r_i=1$ otherwise.}

\subsection{Public scRNA-seq data}

I analyzed a public scRNA-seq dataset available from the Gene Expression Omnibus database under accession number GSE247719 \cite{zhang2025}. There are many public scRNA-seq datasets, and this dataset was selected because UMAP coordinates for each cell were provided and scRNA-seq measurements were performed for multiple time points.

The original dataset contained scRNA-seq data from 14 organs or tissues in male and female mice of three genotypes. In this study, only the data of hepatocytes in liver from male and female wild-type C57BL/6 mice were analyzed. There were five time points: 3, 6, 12, 16, and 23 months. The normal life span of C57BL/6 mice is approximately 2 years.

The number of cells was 81,698 (3 months), 159,805 (6 months), 85,725 (12 months), 47,024 (16 months), and 72,534 (23 months) for \diff{males}, and 93,794 (3 months), 144,678 (6 months), 76,066 (12 months), 84,274 (16 months), and 87,141 (23 months) for \diff{females}. None of the conditions had extremely low cell numbers.

The other parameters for the public scRNA-seq data were as follows: $K=20$, $\varepsilon=10^{-10}$, \diff{$\tau=5$,} $\lambda_1=0.005$, $\lambda_2=0.005$, and $H=\{0.5, 0.8, 0.9\}$. \diff{For males, $r_i=-1$ for $i=3$, and $r_i=1$ otherwise. For females, $r_i=-1$ for $i\in\{3,9\}$, and $r_i=1$ otherwise.}

\section{Results}

\subsection{One-dimensional double-well model}

Figure 2 shows the results for estimating the transition matrix $A$ for the one-dimensional double-well model. The shape of the double-well potential $U(x)$ is shown in Fig.~2A. From the time series, spontaneous noise-induced transitions between the two wells were observed (Fig.~2B). The raw continuous state $x$ was then discretized into $\tilde x$ ranging from 1 to 20 (Fig.~2C). The macroscopic properties of the dynamics were preserved. The histogram of the empirical distribution had two peaks corresponding to the local minima of $U(x)$, as shown in Fig.~2D. By calculating the relative frequencies, the transition matrix $A$ was estimated (Fig.~2E). For most columns, the diagonal element had the highest probability, indicating a tendency to maintain the current state. To check the accuracy of $A$, the stationary distribution $z^*$ was calculated by eigenvalue decomposition of $A$, as shown in Fig.~2F. Its shape was nearly identical to the empirical distribution (Fig.~2D), suggesting the validity of the present approach.

\begin{figure}[t]
\centering
\includegraphics[width=0.8\hsize]{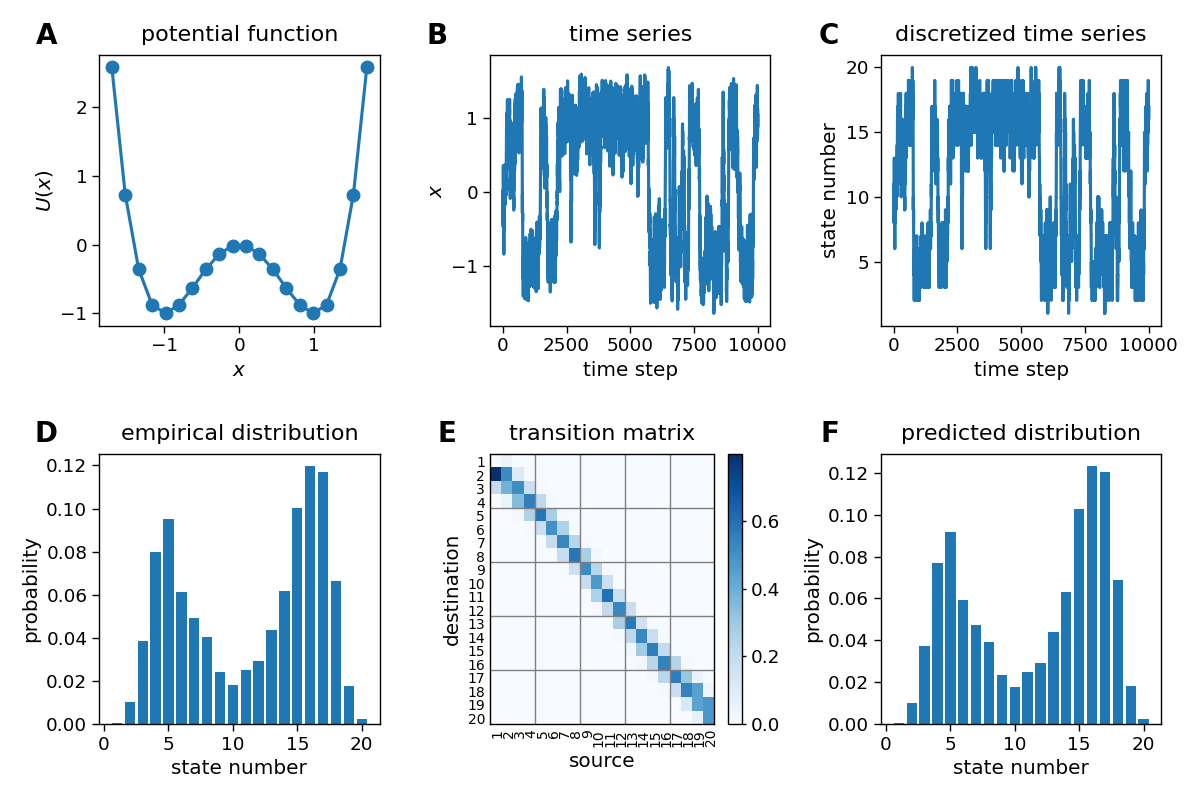}
\caption{Estimation of transition matrix for one-dimensional double-well model. (A) Potential function. Dots indicate the positions of discretized states. (B) Time series. (C) Discretized time series. (D) Empirical distribution. (E) Transition matrix. (F) Predicted distribution.}
\label{fig02}
\end{figure}

Figure 3 shows the effects of smoothing. The length of the time series was shortened to 1,000. In this case, the estimated transition matrix $A$ was slightly noisy without smoothing (Fig.~3A). To resolve this, the Nadaraya-Watson kernel regression, as shown in Eq. (8), was used for smoothing. Figure~3B shows the case of weak smoothing with $\gamma=20$, where $A$ became less noisy. Figure~3C shows the case of strong smoothing with $\gamma=5$, where the concentration along the diagonal elements disappeared. These results suggest that an appropriate level of smoothing may be effective depending on the situation, but excessive smoothing could be harmful.

\begin{figure}[t]
\centering
\includegraphics[width=0.8\hsize]{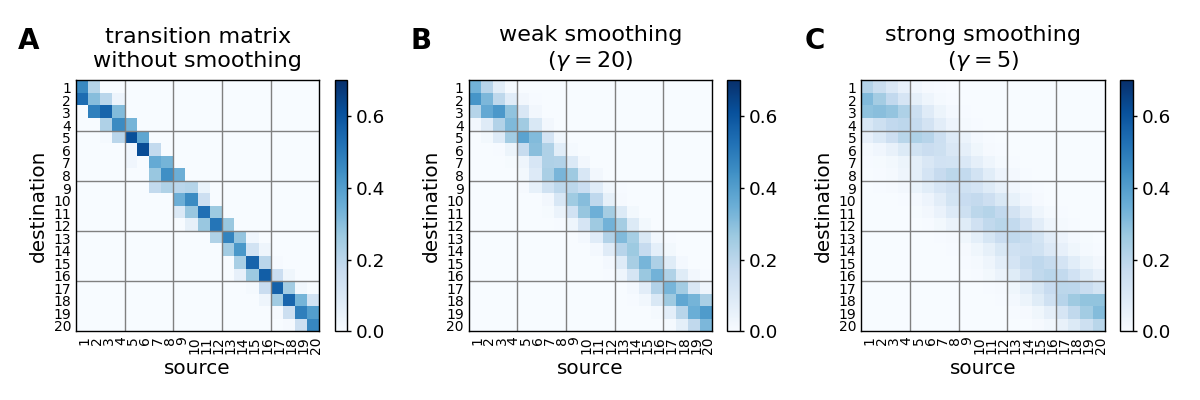}
\caption{Effects of smoothing. (A) Transition matrix without smoothing. (B) Transition matrix with weak smoothing ($\gamma=20$). (C) Transition matrix with strong smoothing ($\gamma=5$).}
\label{fig03}
\end{figure}

Figure 4 shows the results of sparse control for the one-dimensional double-well model. The reward was set to 1 for $\tilde x\leq 10$ and $-1$ for $\tilde x\geq 11$, inducing a bias toward the left well. The magnitudes of the regularization parameters $\lambda_1$ and $\lambda_2$ were varied to see how the result changed. Figure~4A and 4D show the case of $\lambda_1=\lambda_2=0.1$. In this case, only one transition, from state 9 to state 10, was suggested to be suppressed. The reduced flow was added to the self-transition to state 9. The stationary distribution $z^*$ for $A+A'$ is shown in Fig.~4D. We can see that the probability of the left half region increased whereas the probability of the right half region decreased. When the magnitudes of the regularization parameters were decreased to $\lambda_1=\lambda_2=0.05$, two states were suggested to be suppressed (Figs.~4B and 4E). The transition from state 11 to state 12 was added as another intervention candidate. The predicted stationary distribution was further biased toward the left (Fig.~4E). Moreover, when the magnitudes of the regularization parameters were decreased to $\lambda_1=\lambda_2=0.01$, three states were suggested to be suppressed (Figs.~4C and 4F). The transition from state 7 to state 8 was added as a new intervention candidate, and the predicted stationary distribution was almost entirely confined to the left well. These results suggest that the regularization parameters can be used to control the sparsity of $A'$ to a desired level.

When the model’s parameters were changed, states near the center (the saddle point of the potential function) tended to be selected as intervention candidates. This observation may be natural because suppressing the flow near the bottleneck is expected to be most effective.

\begin{figure}[t]
\centering
\includegraphics[width=0.8\hsize]{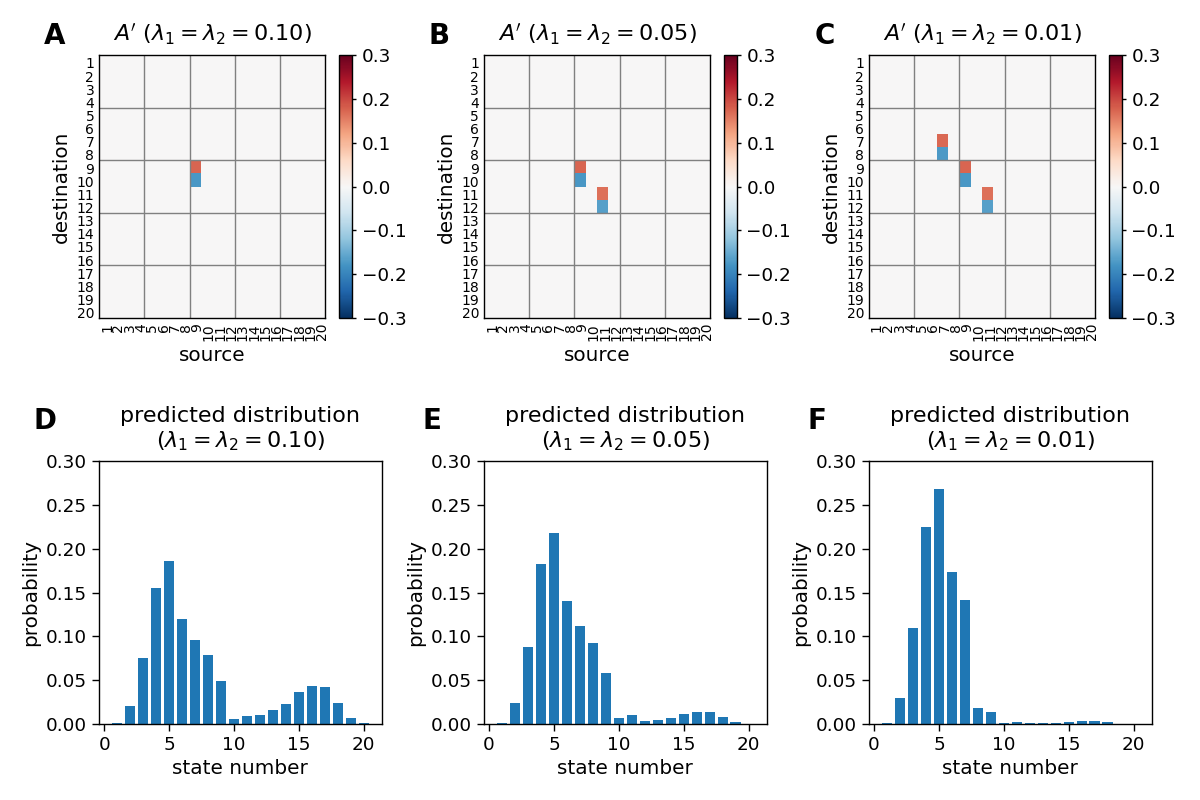}
\caption{Sparse control for one-dimensional double-well model. (A) $A'\ (\lambda_1=\lambda_2=0.1)$. (B) $A'\ (\lambda_1=\lambda_2=0.05)$. (C) $A'\ (\lambda_1=\lambda_2=0.01)$. (D) Predicted distribution $(\lambda_1=\lambda_2=0.1)$. (E) Predicted distribution $(\lambda_1=\lambda_2=0.05)$. (F) Predicted distribution $(\lambda_1=\lambda_2=0.01)$.}
\label{fig04}
\end{figure}

\subsection{Two-dimensional double-well model}

Figure 5 shows the results of MCSC for the two-dimensional double-well model with division per axis. A contour plot of the double-well potential $U(v)$ is shown in Fig.~5A. There were two wells: one at the bottom left and the other at the top right. The trajectory of the stochastic dynamical system showed that lower energy states were visited more frequently (Fig.~5B). From the time series, spontaneous noise-induced transitions between the two wells were observed (Fig.~5C). To discretize the continuous state $v=(x,y)^{\top}$, each axis was divided into 10 regions, and the entire state space was divided into 100 states (Fig.~5D). Although the resolution was not high, we can see that the empirical distribution had two peaks corresponding to the local minima of $U(v)$. By calculating the relative frequencies, the transition matrix $A$ was estimated, and its stationary distribution $z^*$ without control was calculated (Fig.~2E). The predicted distribution was similar to the empirical distribution, supporting the accuracy of $A$. The plot of transition matrix $A$ was omitted because it was difficult to interpret directly.

The reward was set to 1 for \diff{$x\leq 0\land y\leq 0$}, $-1$ for \diff{$x>0\land y>0$}, and 0 otherwise, inducing a bias toward the bottom left well. Many transitions near the saddle point of the potential function were suggested to be suppressed (Fig.~5F). In the predicted stationary distribution with control, the probability of the top right well was successfully reduced, and the probability of the bottom left well was increased. These results suggest that the behavior of MCSC was similar between the one-dimensional double-well model and the two-dimensional double-well model.

When the regularization parameters were changed to reduce the number of intervention sites, the probability of the top right well remained relatively high. This may be because the two-dimensional model has more paths connecting the two wells than the one-dimensional model.

\begin{figure}[t]
\centering
\includegraphics[width=0.8\hsize]{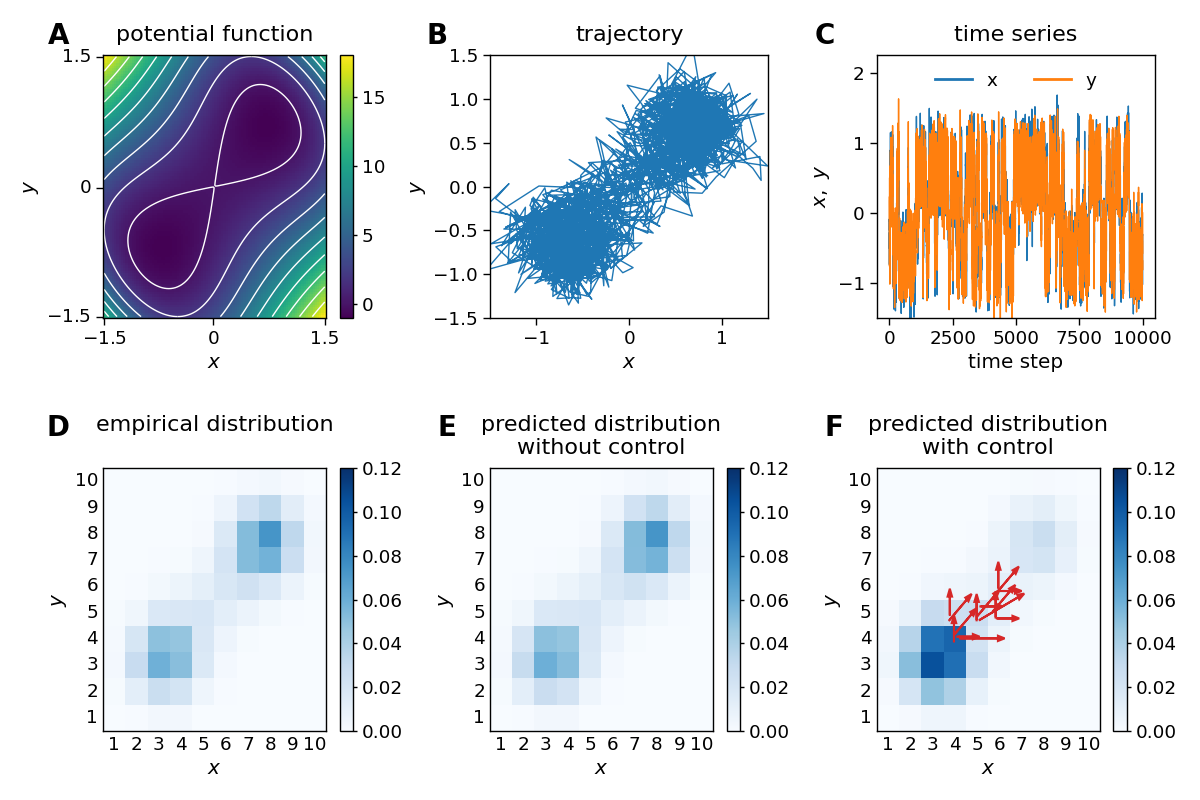}
\caption{MCSC for two-dimensional double-well model with division per axis. (A) Potential function. (B) Trajectory. For clarity, only the first 2000 steps are shown. (C) Time series. (D) Empirical distribution. (E) Predicted distribution without control. (F) Predicted distribution with control. Arrows indicate the suppressed transitions. Color scales for (D--F) indicate the probability.}
\label{fig05}
\end{figure}

Figure 6 shows the results of MCSC for the two-dimensional double-well model with Voronoi partition. Figure 6A shows a scatter plot of the original data points. The data were not exactly the same as those of Fig.~5B due to different random seeds. Figure~6B shows the result of $K$-means partition and the predicted stationary distribution without control. Even though the number of states was only 50, half of that in Fig.~5, the whole state space was efficiently partitioned. This was because the cluster centers were located only in areas with a sufficient number of data points. Suggested interventions were concentrated around the saddle point (Fig.~6C), similar to Fig.~5F. As a control experiment, random selection for representative points was tested, as shown in Figs.~6D and 6E. The area of each Voronoi cell became more heterogeneous, and the suggested interventions became less focused. Given the low computational cost of $K$-means, these results suggest that $K$-means is preferable to random selection. Based on these findings, $K$-means partition was used in the subsequent numerical simulations and real data analysis.

\begin{figure}[t]
\centering
\includegraphics[width=0.8\hsize]{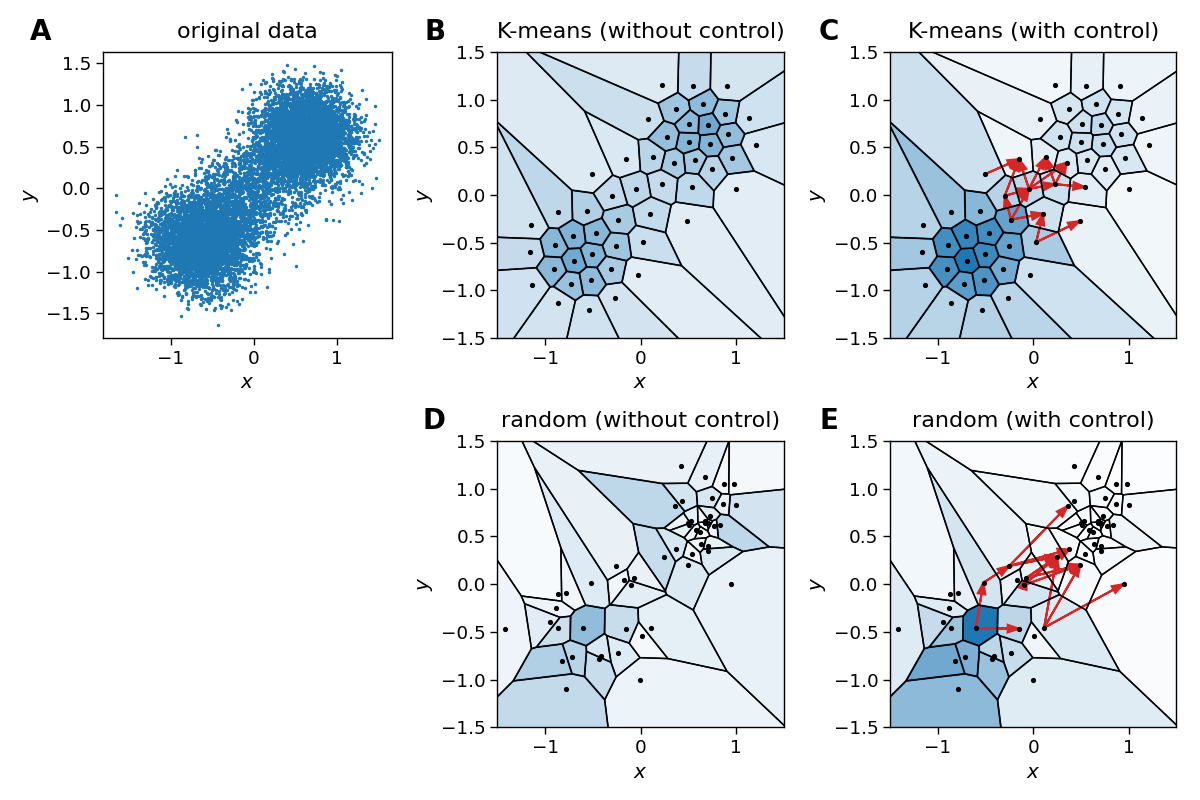}
\caption{MCSC for two-dimensional double-well model with Voronoi partition. (A) Original data points. (B) $K$-means partition and predicted distribution without control. (C) $K$-means partition and predicted distribution with control. (D) Random partition and predicted distribution without control. (E) Random partition and predicted distribution with control. Arrows indicate the suppressed transitions. Color scales for (B--E) indicate the probability.}
\label{fig06}
\end{figure}

\diff{Figure 7 shows the effect of $K$ on MCSC for the two-dimensional double-well model with $K$-means partition. When $K\leq 40$ (Figs.~7A and 7B), the probability of each region was quite heterogeneous, and the reduction of the probability of the top right well was not clear. On the other hand, the results for $K\geq 60$ (Figs.~7C--7E) were similar to each other, indicating that the results converged successfully. On the other hand, the computation time rapidly increased as $K$ increased, as shown in Fig.~7F. These results suggest that sufficiently large but not too large $K$ should be used for MCSC.}

\begin{figure}[t]
\centering
\includegraphics[width=0.8\hsize]{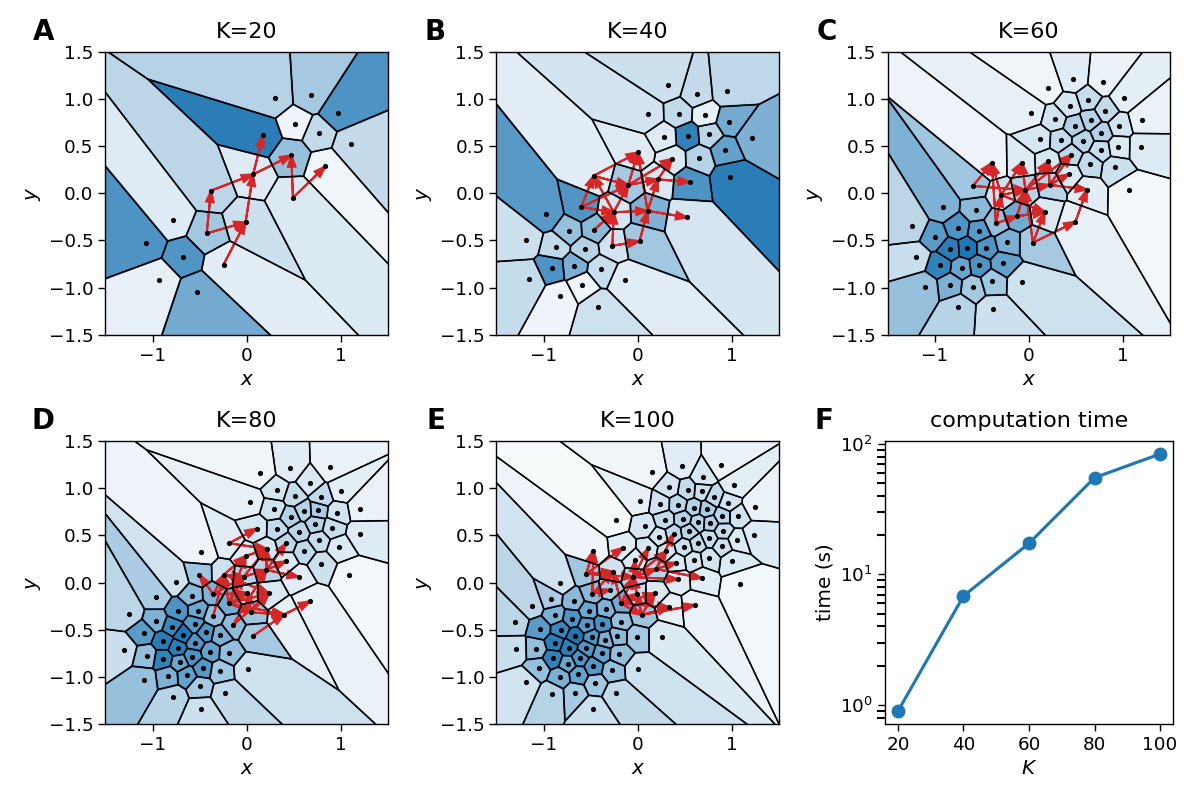}
\caption{\diff{Effect of $K$ on MCSC for two-dimensional double-well model with $K$-means partition. (A) $K=20$. (B) $K=40$. (C) $K=60$. (D) $K=80$. (E) $K=100$. (F) Computation time. Original data points were the same as Fig.~6. Arrows indicate the suppressed transitions. Color scales for (A--E) indicate the probability.}}
\label{fig07}
\end{figure}

\subsection{Two-dimensional branching-flow model}

Figure \diff{8} shows the results of MCSC for the two-dimensional branching-flow model. A heatmap plot of the potential function is shown in Fig.~\diff{8}A. In each trial, the state $v=(x,y)^{\top}$ started at the top center point marked with an asterisk `*’, moved downward while avoiding high energy states, and finally reached one of the four target states. Precisely, the potential function only affected the motion in the $x$ direction, and a constant velocity was assumed in the $y$ direction. Figure~\diff{8}B shows the trajectories for 100 trials, showing that the final destination was determined by stochastic fluctuations. Figure~\diff{8}C shows the $K$-means partition and the predicted stationary distribution without control. Figures~\diff{8}D to \diff{8}G show the control results for each of the four target states. \diff{The reward was set to 1 for the discretized state closest to the target state and $-1$ otherwise.} In all cases, suggested interventions were concentrated on two areas, one for the upstream branching, and the other for the downstream branching. This behavior aligns with intuition, and these results demonstrate that MCSC can automatically produce intuitive and natural outcomes.

\begin{figure}[t]
\centering
\includegraphics[width=0.8\hsize]{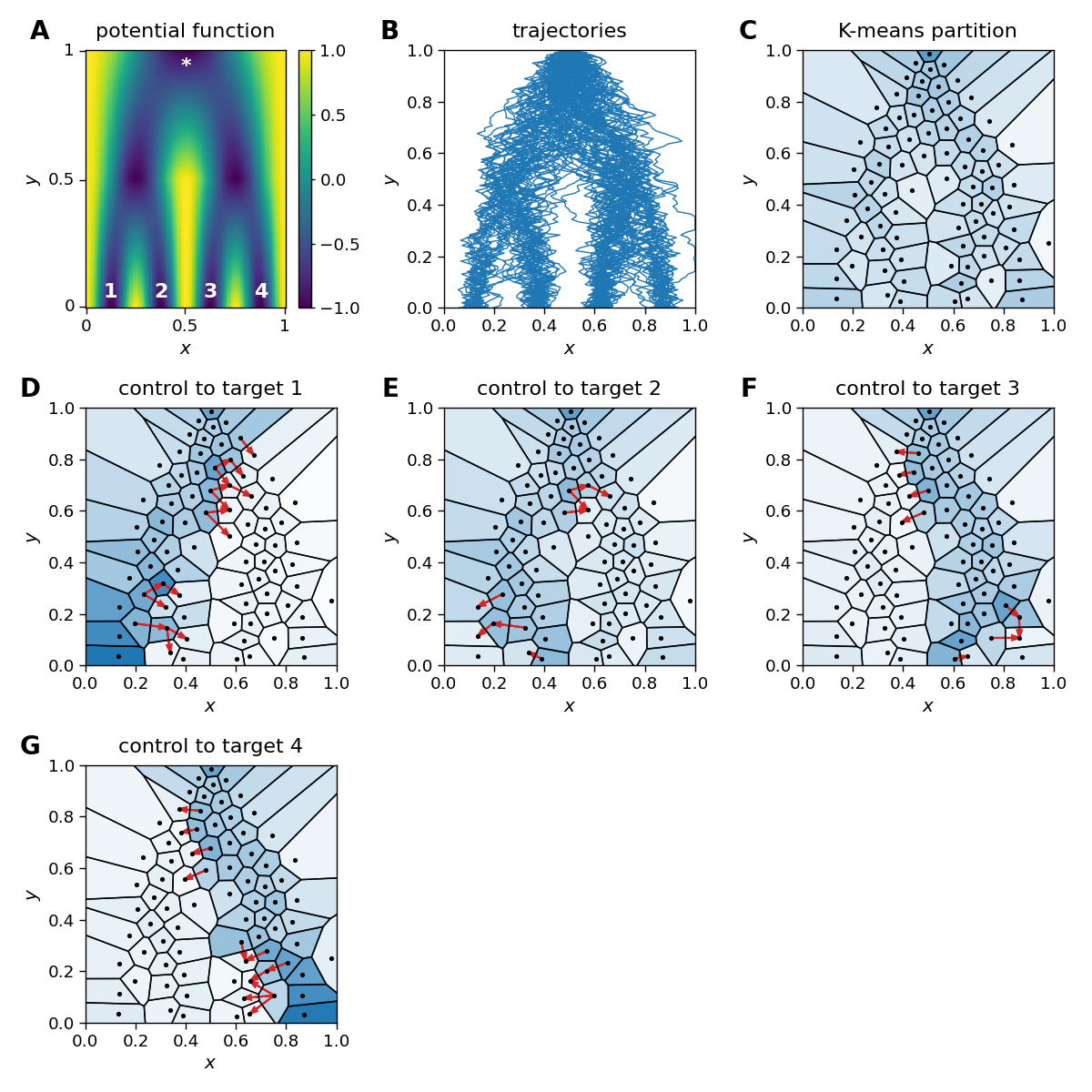}
\caption{MCSC for two-dimensional branching-flow model. (A) Potential function. The asterisk `*’ indicates the start point. Numbers indicate the target states. (B) Trajectories. (C) $K$-means partition and predicted distribution without control. (D-G) Predicted distributions with control to target state 1 (D), target state 2 (E), target state 3 (F), and target state 4 (G). Arrows indicate the suppressed transitions. Color scales for (C--G) indicate the probability for the predicted distributions.}
\label{fig08}
\end{figure}

\subsection{Three-dimensional chaotic attractors}

Figure \diff{9} shows the results of MCSC for the three-dimensional chaotic attractors. Although MCSC primarily targets stochastic dynamical systems, these deterministic models were selected to demonstrate its \diff{applicability} to complex dynamics. Figure~\diff{9}A shows the Lorenz attractor and its $K$-means partition. The strange attractor having two wings was successfully discretized into 50 regions. \diff{Note} that the same color was assigned to multiple clusters due to the difficulty of preparing 50 easily distinguishable colors. For each state $v=(x,y,z)^{\top}$, an integer corresponding to the nearest cluster center was assigned as $\tilde v=\arg\min_k d(v,v^{(k)})$, where $v^{(1)},\ldots,v^{(K)}$ were the cluster centers and $d$ was the Euclidean distance. From the resulting symbolic dynamics, the transition matrix $A$ was calculated, and MCSC was applied. As an example, the reward was set to 1 for \diff{$x>0$} and $-1$ otherwise, attempting to restrict the dynamics to only within the right wing. Figure~\diff{9}B shows the suggested interventions for the Lorenz attractor. We can see that transitions from the right wing to the left wing were to be suppressed. Figure~\diff{9}C shows the Rössler attractor and its $K$-means partition. The reward was set to 1 for \diff{$z<1$} and $-1$ otherwise, attempting to avoid the $z$-directional spikes (intermittent large increases and decreases). Figure~\diff{9}D shows the suggested interventions for the Rössler attractor. We can see that transitions to the onset of $z$-directional spikes were to be suppressed. Interestingly, more upstream interventions were also suggested. A possible explanation is that suppressing outer orbits might guide the flow to inner orbits, contributing to the avoidance of $z$-directional spikes. These results suggest that MCSC can be applied to complex dynamical systems.

\begin{figure}[t]
\centering
\includegraphics[width=0.7\hsize]{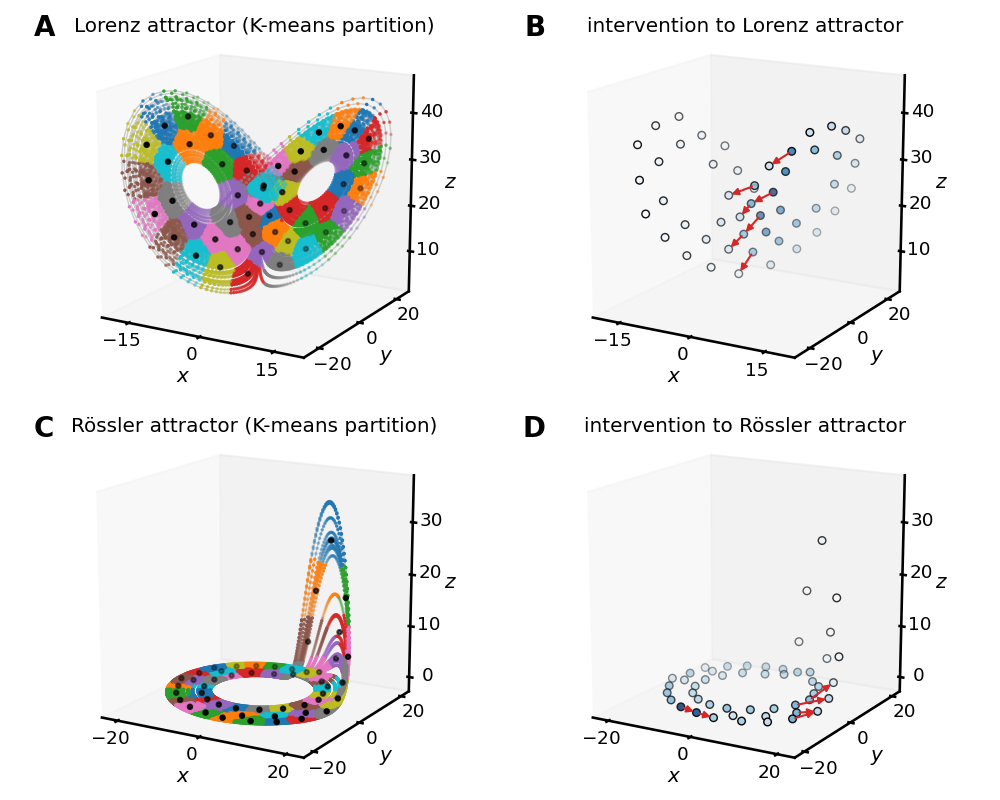}
\caption{MCSC for three-dimensional chaotic attractors. (A) Lorenz attractor and its $K$-means partition. (B) Intervention to Lorenz attractor. (C) Rössler attractor and its $K$-means partition. (D) Intervention to Rössler attractor. Arrows indicate the suppressed transitions. Color scales for (B) and (D) indicate the probability.}
\label{fig09}
\end{figure}

\subsection{Public health checkup data}

Figure \diff{10} shows PCA plots of the public health checkup data obtained from the ninth NDB open data. Males and females were clearly separated (Fig.~\diff{10}A). The data points were arranged by age for each gender (Fig.~\diff{10}B). These results suggest that males and females should be analyzed separately and that it is reasonable to treat these data as time series of virtual individuals.

\begin{figure}[t]
\centering
\includegraphics[width=0.7\hsize]{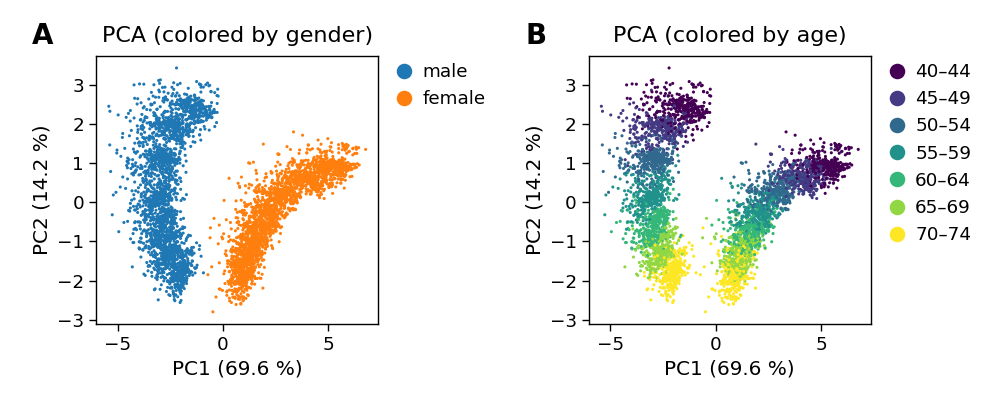}
\caption{PCA plots of public health checkup data. (A) PCA plot colored by gender. (B) PCA plot colored by age.}
\label{fig10}
\end{figure}

Figure \diff{11} shows the results of MCSC for the male data of the public health checkups. The original 12-dimensional state space was divided into 10 regions by using $K$-means. The cluster center of each state is shown in Fig.~\diff{11}A. The states were sorted in advance by the average age of belonging data points. We can see that FPG, HbA1c, and SBP had increasing trends, whereas LDL and ALT had decreasing trends. The other variables showed non-monotonic changes. Figure~\diff{11}B shows the state distribution by age group. State 1 and state 2 were initially dominant in males aged 40--44 years. The dominant states shifted with age, and finally state 9 and state 10 were dominant in males aged 70--74 years. As an example, the reward was set to $-1$ for state 10 and 1 otherwise, attempting to avoid state 10.

\begin{figure}[t]
\centering
\includegraphics[width=0.8\hsize]{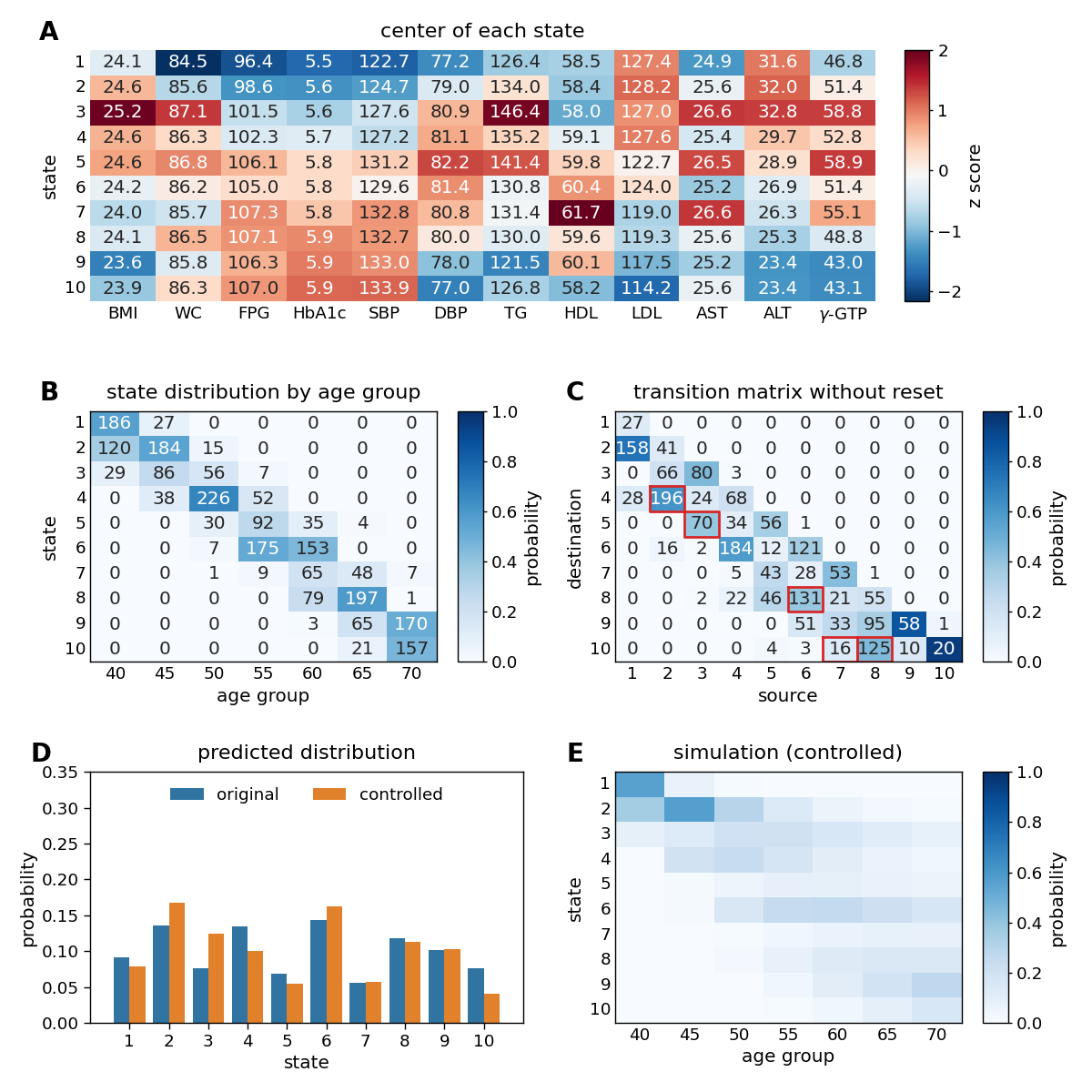}
\caption{MCSC for male data of public health checkups. (A) Center of each state. (B) State distribution by age group. (C) Transition matrix without resetting. Cells with borderlines indicate transitions that were recommended by MCSC to be suppressed. (D) Predicted stationary distributions for original and controlled cases. (E) Simulation of controlled dynamics. Integers in cells in (B) and (C) indicate the number of observations.}
\label{fig11}
\end{figure}

\begin{figure}[t]
\centering
\includegraphics[width=0.8\hsize]{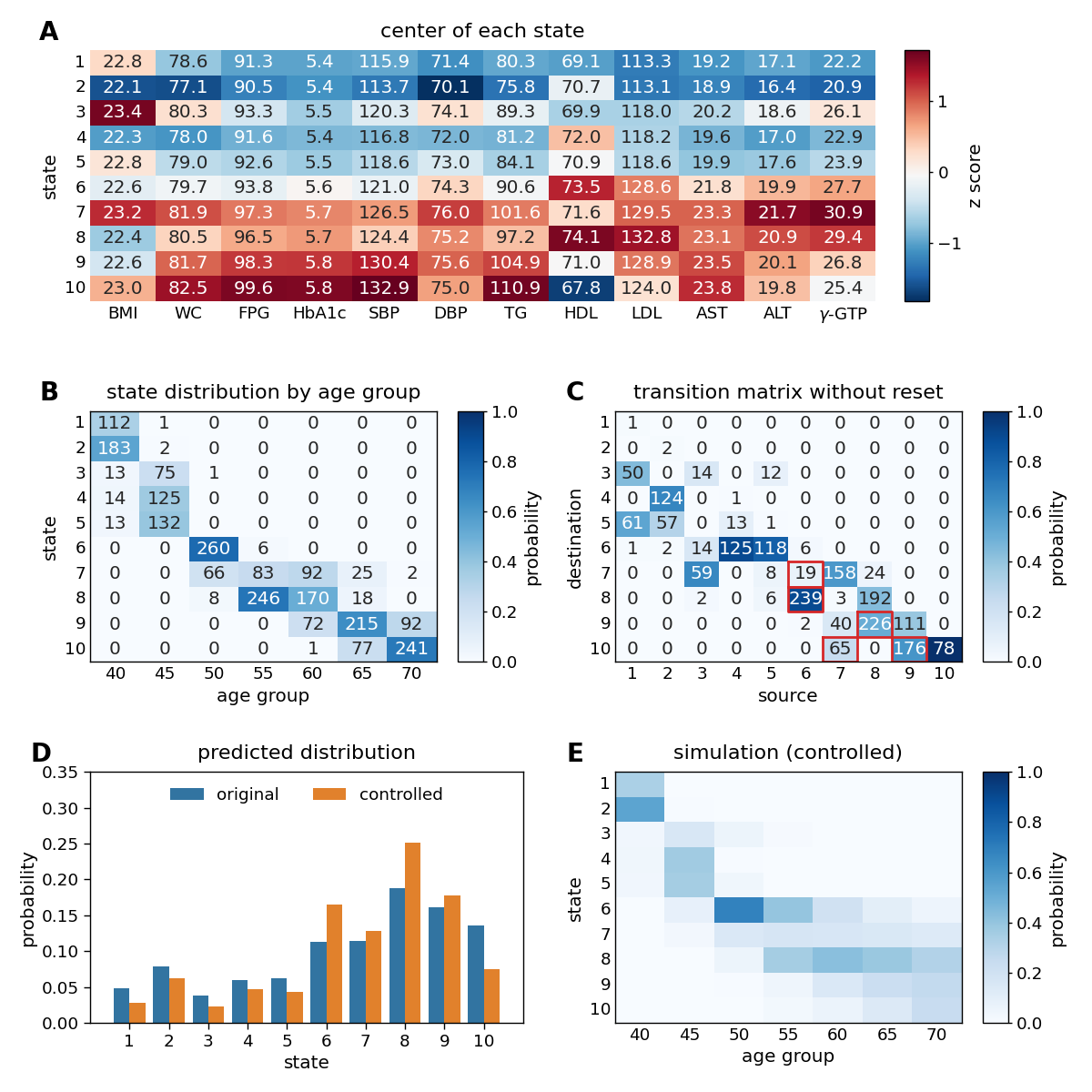}
\caption{MCSC for female data of public health checkups. (A) Center of each state. (B) State distribution by age group. (C) Transition matrix without resetting. Cells with borderlines indicate transitions that were recommended by MCSC to be suppressed. (D) Predicted stationary distributions for original and controlled cases. (E) Simulation of controlled dynamics. Integers in cells in (B) and (C) indicate the number of observations.}
\label{fig12}
\end{figure}

Figure~\diff{11}C shows the transition matrix without resetting. \diff{Note} that MCSC was actually applied to the transition matrix with resetting because otherwise the predicted stationary distribution did not match the empirical distribution. Figure~\diff{11}C is shown only to aid interpretation. The main inflow to state 10 was from state 8, which was suggested to be suppressed. Similarly, the main inflow to state 8 was from state 6, which was also suggested to be suppressed. Interestingly, upstream interventions were also suggested---from state 2 to state 4, and from state 3 to state 5. Because state 2 and state 3 were mainly observed in males aged 40--54 years, these results suggest the effectiveness of early interventions for males. 

Figure \diff{11}D shows the predicted stationary distributions for the original and controlled cases. The probability of state 10 was approximately halved, and the probabilities of state 2 and state 3 were largely increased. Interestingly, the probability of state 4 was also largely decreased, probably due to the upstream intervention. Figure~\diff{11}E shows the result of a simulation of the controlled dynamics. Starting from the initial distribution $z(1)$, the time evolution was calculated as $z(t)=(A+A')^{t-1}z(1)$ $(t=2,\ldots,T)$. The overall flows were dispersed due to interventions.

Figure \diff{12} shows the results of MCSC for the female data of the public health checkups. The data analysis procedures were the same as those in Fig.~\diff{11}. The cluster center of each state is shown in Fig.~\diff{12}A. We can see that WC, FPG, HbA1c, SBP, TG, and AST had increasing trends. DBP, LDL, ALT, and $\gamma$-GTP showed non-monotonic changes peaked at state 7 or state 8. Figure~\diff{12}B shows the state distribution by age group. State 1 and state 2 were initially dominant in females aged 40--44 years. The dominant states shifted with age, and finally state 9 and state 10 were dominant in females aged 70--74 years. As an example, the reward was set to $-1$ for state 10 and 1 otherwise, attempting to avoid state 10.

Figure~\diff{12}C shows the transition matrix without resetting. The main inflow to state 10 was from state 9, which was suggested to be suppressed. Similarly, the main inflow to state 9 was from state 8, which was also suggested to be suppressed. Furthermore, the main inflow to state 8 was from state 6, which was also suggested to be suppressed. In contrast to the male data, upstream interventions were not suggested for females. A possible explanation is that in the case of females, it was almost impossible to stay in states 1 to 5 after 50 years old (Fig.~\diff{12}B), and thus early interventions in states 1 to 5 would have been ineffective. 

Figure \diff{12}D shows the predicted stationary distributions for the original and controlled cases. The probability of state 10 was approximately halved, and the probabilities of state 6 and state 8 were largely increased. The probabilities of states 1 to 5 were also decreased, probably due to resetting. Since this cannot occur in reality, resolving this issue is considered future work. Figure~\diff{12}E shows the result of a simulation of the controlled dynamics. The flows after 50 years old were dispersed due to interventions.

These results demonstrate that MCSC is applicable to the real data of health checkups. In addition, it could suggest a combination of early and late interventions (Fig.~\diff{11}) or only late interventions (Fig.~\diff{12}) depending on the context.

\subsection{Public scRNA-seq data}

Figure \diff{13} shows the UMAP plots of the public scRNA-seq data (GSE247719) \cite{zhang2025}. Only the data of hepatocytes in liver from wild-type mice were used. There were five time points for males (Fig.~\diff{13}A) and females (Fig.~\diff{13}B). The distributions were slightly different at different time points. Although the UMAP coordinates were shared between males and females, their distributions were largely different. Therefore, they were separately analyzed in this study. The state space was divided into 20 regions by using $K$-means for males (Fig.~\diff{13}C) and for females (Fig.~\diff{13}D). In both cases, state 3 was observed only at 23 months. As an example, the reward was designed to prevent the emergence of state 3. Since all analyzed cells were hepatocytes, the unique characteristics of those in state 3 remain unknown. This example was chosen only to demonstrate the effectiveness of MCSC for scRNA-seq data.

\begin{figure}[t]
\centering
\includegraphics[width=0.8\hsize]{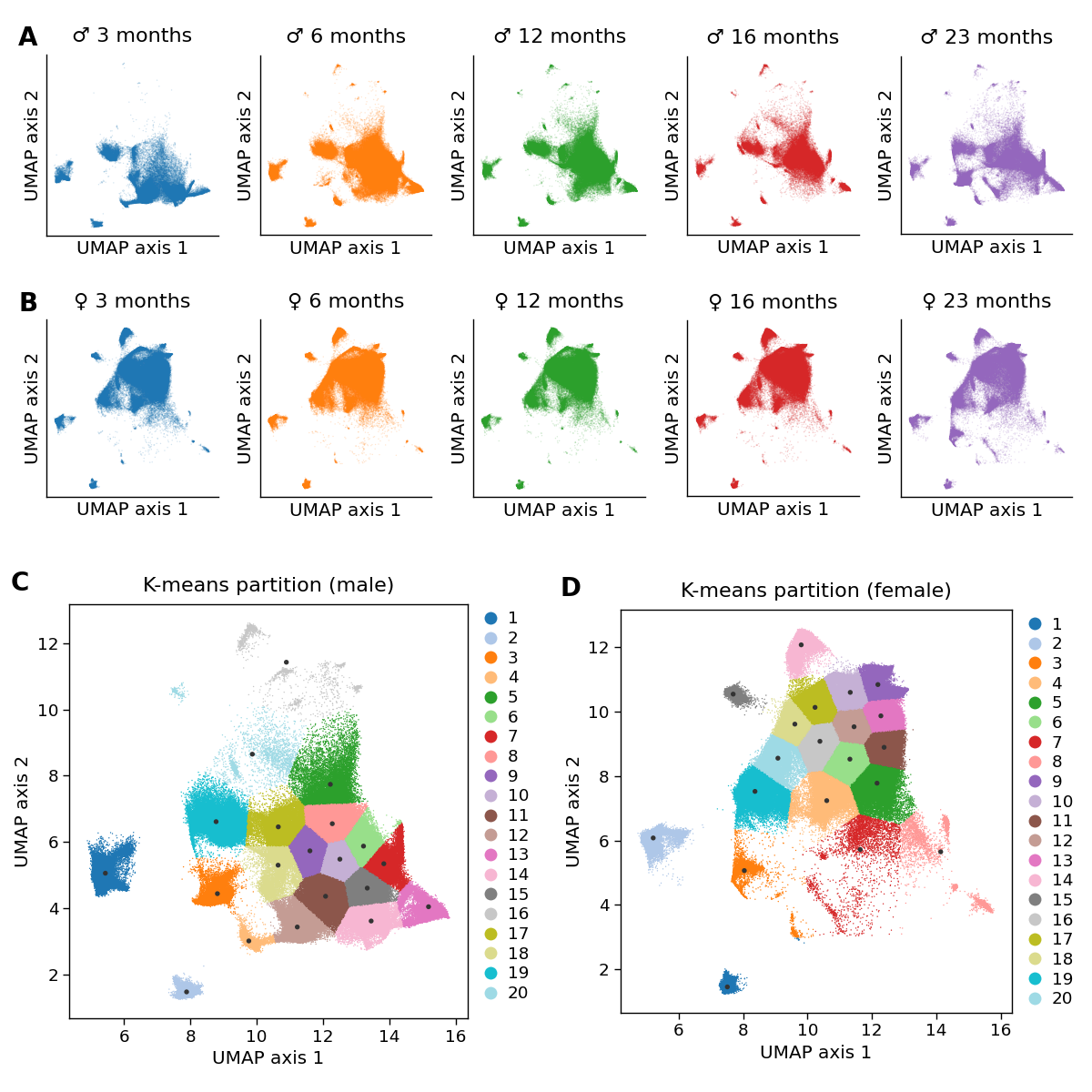}
\caption{UMAP plots of public scRNA-seq data. (A) Male’s UMAP plots at each time point. (B) Female’s UMAP plots at each time point. (C) $K$-means partition of male data. (D) $K$-means partition of female data.}
\label{fig13}
\end{figure}

\begin{figure}[t]
\centering
\includegraphics[width=0.8\hsize]{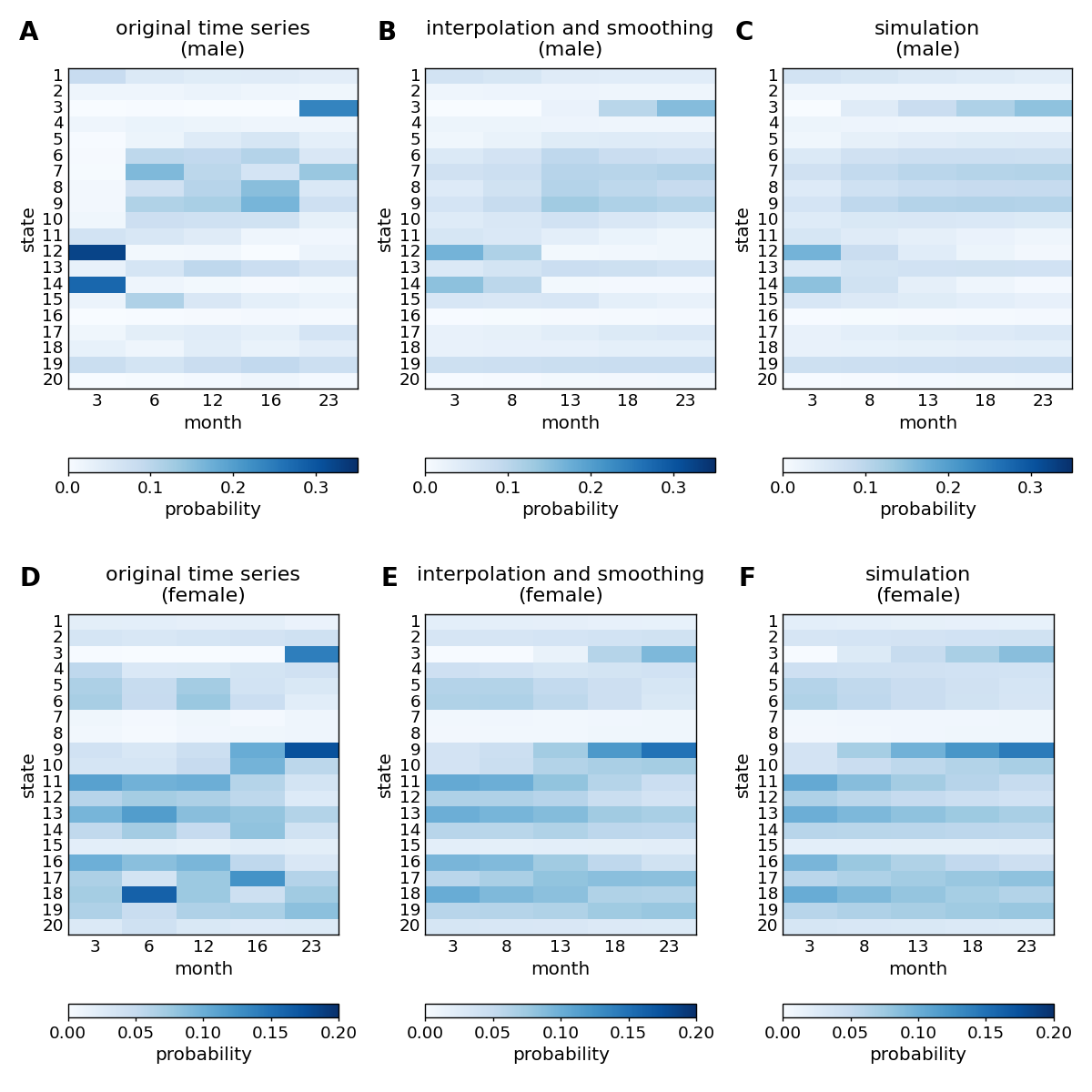}
\caption{Time series of public scRNA-seq data. (A) Original time series for male. (B) Time series with interpolation and smoothing for male. (C) Simulation for male. (D) Original time series for female. (E) Time series with interpolation and smoothing for female. (F) Simulation for female.}
\label{fig14}
\end{figure}

\begin{figure}[t]
\centering
\includegraphics[width=0.8\hsize]{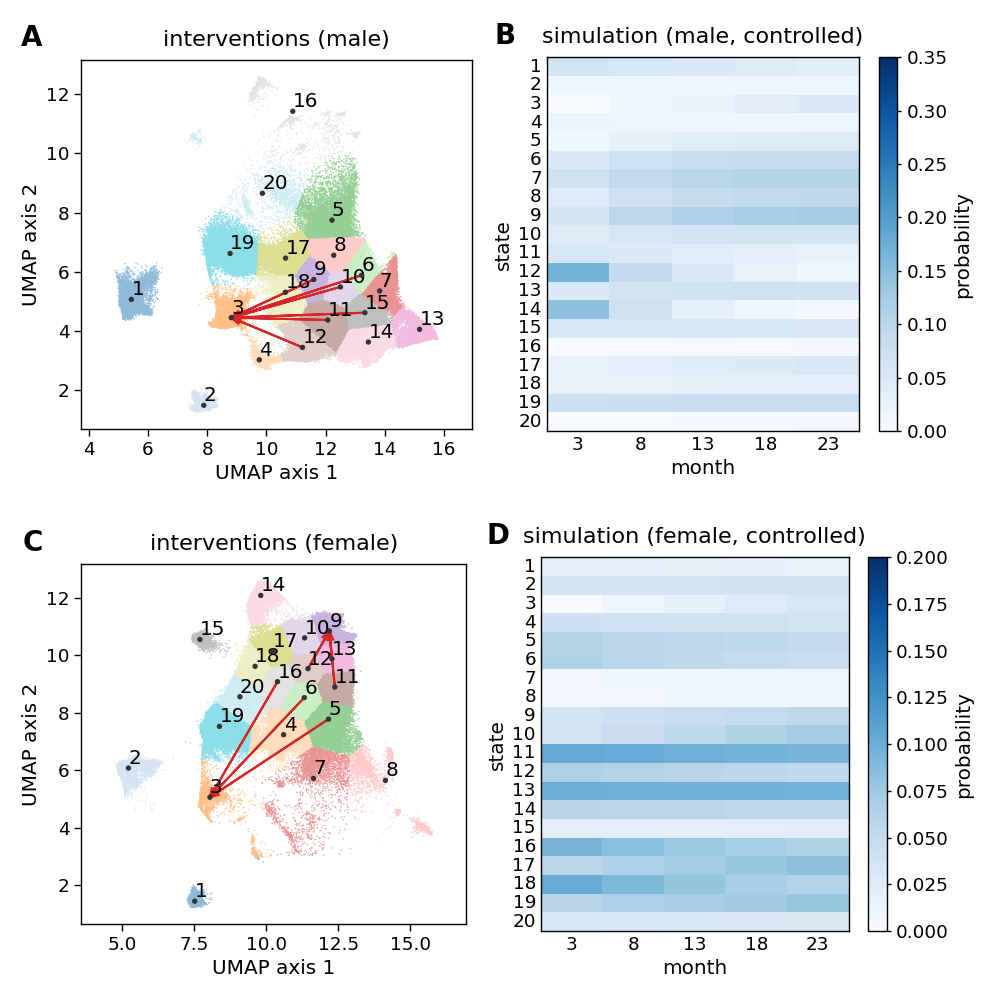}
\caption{Sparse control for public scRNA-seq data. (A) Suggested interventions for male. (B) Simulation of controlled dynamics for male. (C) Suggested interventions for female. (D) Simulation of controlled dynamics for female. Arrows in (A) and (C) indicate the transitions suggested to be suppressed.}
\label{fig15}
\end{figure}

Figure \diff{14} shows the time series of the public scRNA-seq data. The probability distribution $z(t)$ was calculated for each time point. Figures~\diff{14}A and \diff{14}D show the original time series for males and females, respectively. The emergence of state 3 was one of the characteristics of the 23-month time point for both males and females. \diff{Another hallmark was the increased probability of female’s state 9. To suppress it, a negative reward of $-1$ was also assigned to state 9 only for females.} Because the time points were not evenly spaced, they were interpolated and then smoothed (Figs.~\diff{14}B and \diff{14}E). The probability of state 3 was spread to earlier time points. The transition matrix $A$ was estimated using OT, as described in Section 2.3.7. Resetting was not effective for these data and thus omitted. The matrix $A'$ was optimized with respect to a finite horizon $\tau=5$ rather than the stationary distribution. The results of free-run simulations of the uncontrolled dynamics, as shown in Figs.~\diff{14}C and \diff{14}F, resembled Figs.~\diff{14}B and \diff{14}E, respectively. 

Figure \diff{15} shows the results of sparse control for the public scRNA-seq data. Interestingly, long-range transitions were suggested for suppression in both males (Fig.~\diff{15}A) and females (Fig.~\diff{15}C). Remember that the transition matrix was estimated using OT, and short-range transport was preferred. These results suggest that MCSC could suggest intervention candidates that might be difficult for humans to predict. Simulation results showed successful suppression of state 3 emergence in both males (Fig.~\diff{15}B) and females (Fig.~\diff{15}D). \diff{State 9 was also successfully suppressed in females.}

These results demonstrate that MCSC is applicable to real scRNA-seq data using OT-based matching.

\section{Discussion}

The validity of MCSC was demonstrated using numerical simulations of the one-dimensional double-well model (Figs.~2--4), two-dimensional double-well model (Figs.~5--\diff{7}), two-dimensional branching-flow model (Fig.~\diff{8}), and three-dimensional chaotic attractors (Fig.~\diff{9}), as well as real data analyses of the public health checkup data (Figs.~\diff{10}--\diff{12}) and public scRNA-seq data (Figs.~\diff{13}--\diff{15}). In most cases, MCSC was able to automatically suggest interventions that were intuitive and natural to the human eye. For double-well models, interventions near the saddle point were suggested (Figs.~4--\diff{7}). For the branching-flow model, both upstream and downstream interventions were suggested depending on the target state (Fig.~\diff{8}). For the Lorenz attractor, interventions were suggested at the transitions from the right wing to the left wing (Fig.~\diff{9}B). For the Rössler attractor, interventions were suggested around the onset of $z$-directional spikes as well as upstream areas (Fig.~\diff{9}D). For the public health checkup data, both early and late interventions were suggested for males (Fig.~\diff{11}) whereas only late interventions were suggested for females (Fig.~\diff{12}). This difference seems reasonable given the difference in the empirical state distribution by age group for males and females. A possible reason for this sexual difference would be menopause, which generally occurs in women around the age of 50.

Counterintuitive results were obtained only in the scRNA-seq data, where long-range transitions were suggested to be suppressed (Fig.~\diff{15}). This implies that MCSC can suggest potential intervention candidates that would be difficult for humans to predict. However, this study used OT without relays, and it is likely that OT with relays \cite{ling2007, oku2020} can decompose the long-range transitions to multiple short-range transitions.

The idea of state discretization is not at all new. Everyone uses this when plotting histograms. $K$-means is widely used in various data analysis. Binarization of each variable is used in ELA \cite{watanabe2014, ezaki2017, masuda2025}. Symbolic dynamics has long been studied in the field of dynamical systems \cite{hirata2023}. \diff{Piecewise affine systems \cite{imura2010}, finite abstractions \cite{tazaki2008}, and finite automata \cite{kobayashi2012} have long been studied in the field of control theory.} In scRNA-seq studies, it is common to combine multiple cells \cite{baran2019, dann2022, yachimura2024}, collectively referred to as pseudobulk analysis \cite{squair2021}.

In this study, state discretization was used not only to reduce computational cost or improve robustness but also to change the viewpoint. Usually, mathematical \diff{modeling} involves identifying key variables and their interactions. Then, a network model is often considered, where vertices correspond to variables and edges correspond to interactions between variables. This approach is useful for understanding the mechanisms but is data demanding. In particular, biological systems such as the human body are extremely complicated. They are nonlinear, nonstationary, nonautonomous, high-dimensional, hierarchical, and containing delays. Therefore, obtaining quantitative models of biological systems is usually impossible. The data discretization resolves this issue at least partially. An alternative network model is considered, where each vertex corresponds to a point in the state space regardless of the dimension, and each edge represents a transition between discretized states. This phenomenological model can be estimated using relatively little data, even when the original dynamical system is complex. A similar technique is used for RL, where MDP only considers the probability of state transitions given the current state and action, discarding the underlying mechanisms.

MCSC may be related to Koopman mode decomposition (KMD) \cite{susuki2016}, but details have not been elucidated. Both approaches transform nonlinear dynamics into linear models, but there are many differences. KMD considers multiple observation functions, which can be infinitely many, whereas MCSC considers a single observation function $g:\mathbb{R}^{N}\to\{1,\ldots,K\}$. In MCSC, \diff{by} combining the observation $g(x)$ from many trajectories, a probability distribution $z$ within the $(K-1)$-dimensional standard simplex is obtained. KMD basically assumes deterministic dynamical systems, whereas MCSC basically assumes stochastic dynamical systems. The linear operator of KMD may have negative values, whereas the transition matrix $A$ in MCSC must be nonnegative. Further research is needed to clarify their relationships.

\section{Conclusions}

In this study, MCSC was proposed for designing efficient interventions for pre-disease states. Although MCSC is composed of simple and classical methods or theories, such as $K$-means, Markov chains, control theory, and sparse regularization, it has great potential to revolutionize the current healthcare system, contributing to the realization of efficient pre-disease treatment. Its effectiveness was extensively demonstrated using the one- and two-dimensional double-well models, the branching-flow model mimicking the Waddington landscape, the Lorenz attractor, the Rössler attractor, the public health checkup data from the NDB open data, and the public scRNA-seq data with accession number GSE247719. The flexibility and versatility of MCSC have the same root as RL because a Markov chain model can be seen as an actionless MDP model.

Future work has many directions, such as improving the optimization method, incorporating promotion of state transitions, refining the theory, designing controllers that operate in the original continuous space, and applying MCSC to a broader range of real-world datasets. Many ideas and techniques from nonlinear science will help advance these studies to address the problems of an aging society.

\paragraph{Acknowledgments}
The author would like to thank all members of Aihara Moonshot project for valuable discussions. \diff{The author would also like to thank the anonymous reviewers for their careful reading and valuable comments.}

\paragraph{Data availability}
The source codes used in this study \diff{are} available in a GitHub repository \url{https://github.com/okumakito/nolta2025a} under Apache 2.0 license.

\paragraph{Funding}
This research was supported by JST Moonshot R\&D Grant Number JPMJMS2021.

\paragraph{Conflicts of interest}
The author declares no competing interests.

\paragraph{Author contribution}
Makito Oku: Conceptualization, Methodology, Software, Investigation, Formal analysis, Visualization, Writing---Original Draft, and Writing---Review \& Editing.

\paragraph{Artificial intelligence tools}
This study mainly used Google translation for assisting manuscript writing. Microsoft Copilot was also used to check the grammar throughout the manuscript.

\appendix
\section{List of mathematical symbols}
The list of mathematical symbols used in this paper is provided as Table~\ref{tabA1}.

\begin{table}[t]
\caption{List of mathematical symbols}
\label{tabA1}
\begin{minipage}[t]{0.5\textwidth}
\vspace{0pt}
\centering
\begin{tabular}{cp{0.7\textwidth}} \hline
    symbol & description \\ \hline
    $a$ & parameter of Rössler attractor \\
    $a^{(i)}$ & $i$-th column of the transition matrix $A$ \\
    $A$ & transition matrix or system matrix \\
    $A'$ & matrix product $BK'$ \\
    $b$ & parameter of Rössler attractor \\
    $B$ & input matrix \\
    $c$ & regularization parameter; parameter of Rössler attractor \\
    $d$ & distance function \\
    d & symbol for total derivative \\
    $D$ & distance matrix \\
    $f$ & any function \\
    $F$ & transport plan in OT \\
    $g$ & function that maps observation $x$ to discretized state $\tilde x$ \\
    $G$ & objective function of MCSC \\
    $h_i$ & membership functions of fuzzy sets\diff{; element of $H$} \\
    $H$ & set of suppression amount \\
    $i$ & any index \\
    $I$ & identity matrix \\
    $j$ & any index \\
    $k$ & any index \\
    $k_i$ & temporal variables in the Runge-Kutta method \\
    $K$ & number of discretized states \\\
    $K'$ & gain matrix \\
    $l$ & any index \\
    $\diff{L}$ & \diff{dimension of control input} \\
    $m$ & index for individuals \\
    $M$ & number of individuals \\
    $N$ & number of observable variables \\
    $O$ & matrix of ones \\
    $p$ & probability vector \\
    $P$ & probability \\
    $q$ & state transition event \\
    $Q$ & set of state transition events \\ \hline
\end{tabular}
\end{minipage}
\hfill
\begin{minipage}[t]{0.5\textwidth}
\vspace{0pt}
\centering
\begin{tabular}{cp{0.7\textwidth}} \hline
    symbol & description \\ \hline
    $r$ & reward vector \\
    $\mathbb{R}$ & real number \\
    $S_k$ & $k$-th subset in $\mathbb{R}^N$\\
    $t$ & time step \\
    $T$ & number of time points \\
    $u$ & control input \\
    $U$ & potential function \\
    $v$ & multivariate state \\
    $\tilde v$ & discretized \diff{multivariate} state \\
    $w$ & weight \\
    $W$ & Wiener process \\
    $x$ & state; the first coordinate \\
    $x^{(k)}$ & $k$-th representative state \\
    $\tilde x$ & discretized state \\
    $y$ & the second coordinate \\
    $z$ & probability distribution on the discretized states; the third coordinate \\
    $z^*$ & stationary distribution of $z$ \\
    $\alpha$ & damping factor \\
    $\beta$ & parameter of Lorenz attractor \\
    $\Delta^{K-1}$ & ($K-1$)-dimensional standard simplex \\
    $\Delta t$ & time interval \\
    $\varepsilon$ & parameter related to the damping factor \\
    $\gamma$ & parameter of radial basis function \\
    $\eta$ & regularization parameter in weighting \\
    $\lambda_1,\lambda_2$ & regularization parameters in designing a controller \\
    $\pi$ & mathematical constant \\
    $\theta$ & parameter of branching-flow model \\
    $\rho$ & parameter of Lorenz attractor \\
    $\sigma$ & noise intensity; parameter of Lorenz attractor \\
    \diff{$\tau$} & \diff{finite horizon} \\
    $\xi$ & random noise \\ \hline
\end{tabular}
\end{minipage}
\end{table}

\section{Alternative form of Eq.~(\diff{10})}\label{seca2}

If we introduce $\tilde z=(z_1,\ldots,z_{K-1})^{\top}$, $\tilde A\in\mathbb{R}^{(K-1)\times(K-1)}$ with $\tilde A_{ij}=A_{ij}\ (i,j=1,\ldots,K-1)$, $a=(A_{1K},\ldots,A_{K-1,K})^{\top}$, and $e=(1,\ldots,1)^{\top}\in\mathbb{R}^{K-1}$, Eq.~(\diff{10}) can be rewritten as follows:
\begin{equation}
\left[\begin{array}{c}\tilde z(t+1)\\ 1-e^{\top}\tilde z(t+1)\end{array}\right]=
\left[\begin{array}{cc}\tilde A&a\\ e^{\top}-e^{\top}\tilde A&1-e^{\top}a\end{array}\right]
\left[\begin{array}{c}\tilde z(t)\\ 1-e^{\top}\tilde z(t)\end{array}\right].
\end{equation}

This can be reduced to the following equation:
\begin{equation}
\tilde z(t+1)=\tilde A\,\tilde z(t)+a(1-e^{\top}\tilde z(t))=(\tilde A-a\,e^{\top})\tilde z(t)+a.
\end{equation}

By substituting $\tilde z(t),\,\tilde z(t+1)\leftarrow \tilde z^*$, the stationary distribution can be calculated as:
\begin{equation}
\tilde z^*=(I-\tilde A+a\,e^{\top})^{-1}a,
\end{equation}
where $I$ is a $(K-1)$-dimensional identity matrix. By using $\tilde z^*$, the reduced equation can be written as follows:
\begin{equation}
\tilde z(t+1)-\tilde z^* = (\tilde A-a\,e^{\top})(\tilde z(t)-\tilde z^*).
\end{equation}

By introducing $\tilde z'=\tilde z-\tilde z^*\in\mathbb{R}^{K-1}$ and $\tilde A'=\tilde A-a\,e^{\top}\in\mathbb{R}^{(K-1)\times(K-1)}$, an alternative form of Eq.~(\diff{10}) with its origin being the unique equilibrium point can be obtained as follows:
\begin{equation}
\tilde z'(t+1)=\tilde A'\tilde z'(t).
\end{equation}

\bibliographystyle{IEEEtran}
\bibliography{references}

\end{document}